\documentclass{article}
\usepackage[utf8]{inputenc}
\usepackage[english]{babel}
 
\usepackage{amsmath}
\usepackage{amsthm}
\usepackage{tikz}
\usetikzlibrary{svg.path}
\usetikzlibrary{shapes,arrows,chains}
\usepackage{amsfonts}
\usepackage{amssymb}
\usepackage{graphicx}
\usepackage[a4paper,left=3cm,right=3cm]{geometry}
\usepackage{xcolor}
\usepackage{bbm}
\usepackage{csquotes}
\usepackage{latexsym}
\usepackage{float}
\usepackage{hyperref}
\usepackage[backend=biber, style=alphabetic, sorting=nyt, maxbibnames=99]{biblatex}
\newtheorem{thm}{Theorem}[section]
\newtheorem*{thm*}{Theorem}
\newtheorem{lem}[thm]{Lemma}

\newtheorem{prop}[thm]{Proposition}
\newtheorem{coro}[thm]{Corollary}
\newtheorem{conj}[thm]{Conjecture}

\theoremstyle{definition}
\newtheorem{defn}[thm]{Definition}

\newtheorem{rem}[thm]{Remark}
\newtheorem{exa}[thm]{Example}
\newtheorem{exas}[thm]{Examples}

\newtheorem{fact}[thm]{Fact}
\newtheorem{nota}[thm]{Notation}
\newcommand{\rk}{\textnormal{rk}}
\newcommand{\trdeg}{\textnormal{td}}

\newcommand{\ecl}{\textnormal{ecl}}
\newcommand{\cl}{\textnormal{cl}}

\newcommand{\ldim}{\textnormal{ldim}}

\renewcommand{\exp}{\textnormal{exp}}
\newcommand{\gen}[1]{\langle #1 \rangle}

\newcommand{\im}{\textnormal{im}}

\newcommand{\qftp}{\textnormal{qftp}}

\renewcommand{\epsilon}{\varepsilon}
\renewcommand{\phi}{\varphi}

\newcommand{\pcl}{\textnormal{pcl}}

\newcommand{\pdim}{K\textnormal{-td}}
\newcommand{\edim}{\textnormal{etd}}
\newcommand{\Loc}{K\textnormal{-pLoc}}
\newcommand{\treqcl}{\triangleleft_{\textnormal{cl}}}
\newcommand{\treq}{\triangleleft}
\newcommand{\hull}[1]{\lceil #1 \rceil}

\newcommand{\alg}{\textnormal{alg}}
\newcommand{\tr}{\textnormal{tr}}

\newcommand{\affloc}{\textnormal{AffLoc}}
\newcommand{\algloc}{\textnormal{Loc}}

\def\Ind#1#2{#1\setbox0=\hbox{$#1x$}\kern\wd0\hbox to 0pt{\hss$#1\mid$\hss}
\lower.9\ht0\hbox to 0pt{\hss$#1\smile$\hss}\kern\wd0}

\def\Notind#1#2{#1\setbox0=\hbox{$#1x$}\kern\wd0\hbox to 0pt{\mathchardef
\nn="3236\hss$#1\nn$\kern1.4\wd0\hss}\hbox to 0pt{\hss$#1\mid$\hss}\lower.9\ht0
\hbox to 0pt{\hss$#1\smile$\hss}\kern\wd0}

\title{Quasiminimality of Complex Powers}
\author{Francesco Gallinaro and Jonathan Kirby}

\begin{document}

\maketitle
\let\thefootnote\relax\footnotetext{Both authors were supported by EPSRC grant EP/S017313/1. The first author was partially supported by a London Mathematical Society Early Career Fellowship and by the program GeoMod ANR-19-CE40-0022-01 (ANR-DFG).


}

\begin{abstract}
The complex field, equipped with the multivalued functions of raising to each complex power, is quasiminimal, proving a conjecture of Zilber and providing evidence towards his stronger conjecture that the complex exponential field is quasiminimal. 
\end{abstract}

\tableofcontents

\section{Introduction}

Over 25 years ago, Zilber stated his Quasiminimality Conjecture for complex exponentiation:
\begin{conj}[\cite{Zil97}]\label{qm conj}
	The complex field with the exponential function, $\mathbb{C}_\exp:=  \langle \mathbb{C};+,\cdot,\exp\rangle $, is quasiminimal: every subset of $\mathbb C$ which is definable in this structure is countable or co-countable.
\end{conj}

Definable means in the sense of first-order logic: a formula $\phi(z)$ is built from variables, constants for each complex number (so definable ``with parameters''), the operations $+,\cdot,\exp$, and equality $=$ using the usual Boolean operations and quantifiers. We assume there is only one free variable, $z$, and then $\phi(z)$ defines the set of complex numbers $c$ for which $\phi(c)$ is true. 
If we do not use quantifiers, then a definable subset is just a Boolean combination of zero sets of complex exponential polynomials. Such subsets are easily seen to be countable or co-countable. 
In the complex field $\mathbb{C}_{\textnormal{field}}$, without exponentiation, a subset of $\mathbb{C}$ defined without quantifiers is just a Boolean combination of zero sets of polynomials, so is finite or cofinite. By the Tarski--Chevalley quantifier elimination theorem \cite{Tar51}, using quantifiers gives no new definable subsets in this case, so $\mathbb{C}_{\textnormal{field}}$ is \emph{minimal}: every definable subset is finite or cofinite. Another consequence of quantifier elimination is that the definable (or even interpretable) sets in any number of variables correspond to complex algebraic varieties.

With exponentiation this is far from true. The subset $\mathbb{Z}$ of integers is definable in $\mathbb{C}_\exp$, and so the whole arithmetic hierarchy of subsets of $\mathbb{Z}$ is definable, and their complexity increases with each additional quantifier alternation. However, all these subsets are of course countable. It is not known how complicated the subsets of $\mathbb{C}$ can get with increasing quantifier alternations.
One possibility mentioned in \autocite[p.791]{Mar06} is that the real field $\mathbb{R}$ is a definable subset. If so, with $\mathbb{R}$ and $\mathbb{Z}$ together, every subset of $\mathbb{C}$ or $\mathbb{R}$ in the projective hierarchy of descriptive set theory is definable \cite[Exercise 37.6]{Kec}, in particular all continuous functions $\mathbb{C} \to \mathbb{C}$, so definable sets would generally have nothing to do with the exponential.
If Conjecture~\ref{qm conj} is true, the picture is very different, and the definable sets should have a geometric nature, from complex analytic geometry, much closer to algebraic geometry but also with connections to the transcendence theory and the diophantine geometry related to the complex exponential.

We note that in the real case, the definable sets of the real field are the \emph{semi-algebraic sets} and the structure is so-called \emph{o-minimal}. By Wilkie's very influential theorem~\cite{Wil96}, the real exponential field $\mathbb{R}_\exp$ is also o-minimal, and so the definable sets there are geometric in nature.

\medskip
The Quasiminimality Conjecture has sparked a lot of mathematical activity. For example, Zilber's part of the Zilber--Pink conjecture of diophantine geometry, and the related work on functional transcendence around the Ax--Schanuel theorem, came out of his early work towards his conjecture \cite{Zil02}. 
The Pila--Wilkie theorem for counting integer or rational points is now a major tool in diophantine geometry which also came, at least partly, from work towards the conjecture. Indeed, Wilkie \cite{Wil23} wrote:
\begin{quotation}
\ldots my motivation for studying integer points in o-minimally definable sets was, apart from the fun of it, completely motivated by Boris’ [Zilber's] quasiminimality problem.
\end{quotation}
In a third, more model-theoretic direction, quasiminimality is not a property traditionally studied in model theory because it is not an invariant of the finitary first-order theory of a structure. The importance of this conjecture has led to a resurgence of interest in the use of infinitary logics and abstract elementary classes in model theory, as for example in the books of Baldwin \cite{Bal} and Marker \cite{Mar16}.

\medskip

While in this paper we do not prove Conjecture~\ref{qm conj}, we do prove a strong result towards it, the Quasiminimality Conjecture for complex powers.

\begin{thm}\label{first}
For $\lambda \in \mathbb C$, let $\Gamma_\lambda = \{(\exp(z), \exp(\lambda z)) \mid {z \in \mathbb{C}}\}$ denote the graph of the multivalued map $w \mapsto w^\lambda$. Then the structure
$\mathbb{C}_{\mathbb{C}\text{-powers}} = \langle \mathbb{C};+,\cdot,(\Gamma_\lambda)_{\lambda \in \mathbb{C}}\rangle $ of the complex field equipped with all complex powers is quasiminimal.
\end{thm}

Complex powers are obviously definable in $\mathbb{C}_\exp$, so Theorem~\ref{first} would be a consequence of Conjecture~\ref{qm conj}, and is the most significant result yet proved towards it. We give some discussion of the gap between Theorem~\ref{first} and Conjecture~\ref{qm conj} after the statement of Theorem~\ref{second}.

For the structure $\mathbb{C}_{\mathbb{C}\text{-powers}}$ itself, a consequence of Theorem~\ref{first} is that $\mathbb{R}$ is not definable, and so neither is the projective hierarchy of its subsets. So at least on a generic level, the definable sets are geometric in nature. There is no obvious definition of $\mathbb{Z}$ as a ring in $\mathbb{C}_{\mathbb{C}\text{-powers}}$. Indeed, it follows from \cite{Zil11} that if we restrict to powers from certain subfields $K$, the theory is superstable so $\mathbb{Z}$ as a ring is not interpretable and all the definable sets are geometric in nature. Assuming Schanuel's conjecture, this conclusion holds for all complex powers. However, we cannot prove it unconditionally.

\subsection{Towards the Quasiminimality Conjecture}

Boxall made progress towards Conjecture \ref{qm conj} by showing in ~\cite{Box20} that certain existential formulas in the language of exponential rings must define countable or co-countable sets in $\mathbb{C}_\exp$.

Wilkie has an approach to proving Conjecture~\ref{qm conj} via analytic continuation of definable holomorphic functions, discussed in \cite{Wil13, Wil23}. In 2008 he gave talks  \cite{Wil08c} announcing a proof of the quasiminimality of $\mathbb{C}$ with the power $i$ only, that is, of $\langle \mathbb{C};+,\cdot,\Gamma_i\rangle$. However, the method relied heavily on the fact that $i^2 = -1$, and did not extend to other powers, and the proof has not yet appeared. 

Zilber's own attempts to prove his conjecture led to a construction of a quasiminimal exponential field \cite{Zil05}, now known as $\mathbb{B}_\exp$, via the Hrushovski--Fra\"iss\'e \emph{amalgamation-with-predimension} method \cite{Hru}, which produces a countable structure, and then by his own variant \cite{Zil05b} (see also \cite{Kir10b,BHHKK}) of Shelah's \emph{excellence} method \cite{She}, which extends the amalgamation to uncountable cardinalities, in particular to the continuum-sized model $\mathbb{B}_\exp$. Zilber then conjectured that $\mathbb{C}_\exp \cong \mathbb{B}_\exp$, which implies Conjecture~\ref{qm conj}, but is much stronger, since it incorporates Schanuel's Conjecture of transcendental number theory. In fact, it is equivalent to that conjecture together with another conjectural property for $\mathbb{C}_\exp$, called \emph{Strong Exponential-Algebraic Closedness} (SEAC) which asserts that certain systems of exponential polynomial equations in many variables should have complex solutions, but also includes a condition that the solutions should have large enough transcendence degree.

Bays and the second author \cite{BK18} modified Zilber's construction and were able to remove the transcendence conjectures from this path to Conjecture~\ref{qm conj}. They were able to show that Conjecture~\ref{qm conj} follows from \emph{Exponential-Algebraic Closedness} (EAC) for $\mathbb{C}_\exp$, which is like SEAC but without the transcendence requirement.

The proof of Theorem~\ref{first} in this paper follows this Bays--Kirby strategy for quasiminimality, adapted for complex powers rather than complex exponentiation.

\subsection{Exponential sums with a field of exponents}

We actually prove a slightly stronger form of Theorem~\ref{first}: quasiminimality for the complex field in Zilber's exponential sums language, first introduced in \cite{Zil03}. For a subfield $K \subseteq \mathbb{C}$, we write $\mathbb{C}^K$ for the complex numbers as a \emph{$K$-powered field}, that is, as the 2-sorted structure
$$\mathbb{C}_{K\text{-VS}} \stackrel{\exp}{\longrightarrow} \mathbb{C}_{\textnormal{field}}$$
where the image sort is $\mathbb{C}$ equipped with the field structure, the covering sort is $\mathbb{C}$ equipped only with its structure as a $K$-vector space, and the covering map is the usual complex exponentiation.

\begin{thm}\label{second}
The structure $\mathbb{C}^\mathbb{C}$ is quasiminimal. That is, any definable subset of either sort (in one free variable) is countable or its complement in that sort is countable.
Equivalently, for any countable subfield $K \subseteq \mathbb{C}$, the structure $\mathbb{C}^K$ is quasiminimal.
\end{thm}
The complex power functions $\Gamma_\lambda$ are definable in $\mathbb{C}^\mathbb{C}$ so Theorem~\ref{first} follows immediately from Theorem~\ref{second}. 

In analytic number theory, an exponential sum is an expression of the form $\sum_{j=1}^r a_j \exp(2 \pi i x_j)$ where the $a_j$ are real coefficients and the $x_j$ are real or complex numbers or variables. Generally, one looks to get bounds on such expressions. A consequence of Theorem~\ref{second} is that in $\mathbb{C}^\mathbb{C}$ we cannot define absolute values or indeed the set of reals or the order relation. 
What we can express are complex exponential sums \emph{equations}.  
For example, if $z_1,\dots,z_n$ are variables in the $K$-vector space sort, $\alpha$ is an element of that sort, and $\lambda_1,\dots,\lambda_n \in K$, then $\exp\left(\alpha+ \sum_{i=1}^n \lambda_i z_i\right)$ is a term in the field sort. Writing $w_i=\exp(z_i)$ and $a = \exp(\alpha)$, we can informally write this as $a \prod_{i=1}^n w_i^{\lambda_i}$ to get a ``monomial'' with complex exponents. Applying addition in the field sort, we can then get ``polynomials'' with complex exponents, and these terms, treated properly with the variables $z_i$ from the covering sort, are the exponential sums whose zero sets are generalisations of complex algebraic varieties and which are the basic definable sets in $\mathbb{C}^K$.

However, in $\mathbb{C}^K$ we cannot iterate exponentiation, and in general there is no way to recover the embedding of $K$ as a subfield of the field sort $\mathbb{C}_{\textnormal{field}}$ from the structure $\mathbb{C}^K$. Likewise, there is no way to recover the field structure on the covering sort, let alone identify it with the field sort. If we could do that, we could iterate the exponential map and $\mathbb{C}^\mathbb{C}$ would be bi-interpretable with $\mathbb{C}_\exp$, so Theorem~\ref{second} would actually prove Conjecture~\ref{qm conj}.

If we take $K = \mathbb{Q}$, the structure $\mathbb{C}^\mathbb{Q}$ was axiomatised and shown to be quasiminimal in \cite{Zil06, BZ11}. Since rational powers are algebraic in nature, there is no analytic content to this structure. 
If we now take $K = \mathbb{Q}(\lambda)$ for some $\lambda \in \mathbb{C}$, the structure $\mathbb{C}^K$ depends on $\lambda$. In topological terms, there is a clear difference between real and non-real $\lambda$. For $\lambda \in \mathbb{R}$, there is a branch of $w \mapsto w^\lambda$ which fixes the positive real line setwise, but not for non-real $\lambda$. One might wonder if this distinction shows up in the algebra of powers: for example, if $t$ is real and transcendental, are the structures $\mathbb{C}^{\mathbb{Q}(t)}$ and $\mathbb{C}^{\mathbb{Q}(it)}$ distinguishable, or isomorphic? It turns out that not all transcendental powers are isomorphic, and we give an example in \ref{taugeneric}. However, if Conjecture~\ref{qm conj} is true then all but countably many complex powers should give isomorphic powered fields, and indeed we are able to prove this. 

Given a countable field $K$, we construct a $K$-powered field $\mathbb{E}^K$ of cardinality continuum, analogous to Zilber's $\mathbb{B}_\exp$.
We prove

\begin{thm}\label{third}
    Let $K$ be a countable field of characteristic 0. Then up to isomorphism, there is exactly one $K$-powered field $\mathbb{E}^K$ of cardinality continuum which:
\begin{itemize}
\item[(i)] has cyclic kernel,
\item[(ii)] satisfies the Schanuel property,
\item[(iii)] is $K$-powers closed, and 
\item[(iv)] has the countable closure property.
\end{itemize}
Furthermore, it is quasiminimal.
\end{thm}

In the case $K=\mathbb{Q}(\lambda)$ with $\lambda \in \mathbb{C}$ transcendental, we say that $\lambda$ is a \emph{generic power} if $\mathbb{C}^K \cong \mathbb{E}^K$, so Theorem~\ref{third} implies that generic powers give rise to isomorphic powered fields. We are able to prove:

\begin{thm}\label{fourth}
If $\lambda \in \mathbb C$ is exponentially transcendental then $\lambda$ is a generic power.
\end{thm}

There are only countably many complex numbers which are not exponentially transcendental (that is, which are exponentially algebraic), so this proves the promised corollary of Conjecture~\ref{qm conj}.

We briefly explain the terms used in the statement of Theorem~\ref{third} and how they are proved in the complex case. \emph{Cyclic kernel} just refers to the fact that the kernel of the exponential map is an infinite cyclic group. The \emph{Schanuel property} is a form of Schanuel's conjecture appropriate for powers from a field $K$. It was proved for exponentially transcendental $\lambda$ in \cite{BKW}. There is a natural pregeometry on a $K$-powered field, analogous to relative algebraic closure on a field, which we call \emph{$K$-powers closure}. The \emph{countable closure property} (CCP) asserts that the $K$-powers closure of a countable set is countable. It holds on $\mathbb{C}^K$ because the $K$-powers closure of a set $A$ is contained in the exponential algebraic closure of $A \cup K$, and the exponential algebraic closure has the CCP for topological reasons.

On the other hand the \emph{$K$-powers closed} property does not refer to this pregeometry, but is rather an existential closedness condition: every system of $K$-powers (or exponential sums)  equations which can have solutions in a ``reasonable'' extension of $\mathbb{E}^K$ already has solutions inside $\mathbb{E}^K$. So this is the analogue of a field being (absolutely) algebraically closed, and of the EAC property for exponential fields mentioned above.
Zilber \cite{Zil02} was able to prove this property for $\mathbb{C}^K$ in the case where $K$ was a subfield of $\mathbb{R}$, and then in the unpublished \cite{Zil11} he proved Theorem~\ref{fourth} in the case of real $\lambda$. The main breakthrough which allows us to prove Theorem~\ref{fourth} without this restriction to real powers is:
\begin{fact}[{\autocite[Corollary 8.10]{Gal22}}]\label{C-powers closed fact}
The complex field with complex powers, $\mathbb{C}^\mathbb{C}$, is powers-closed.
\end{fact}

This fact is also essential for the proofs of Theorems~\ref{first} and~\ref{second}. For these we use the Bays--Kirby variation of the Shelah--Zilber excellence method from~\cite{BK18} to construct quasiminimal $K$-powered fields $\mathbb{E}^{K,\tr}(D)$ over a base $K$-powered field $D$, such that $D$ is in fact relatively $K$-powers closed in $\mathbb{E}^{K, \tr}(D)$. Given a countable subfield $K \subseteq \mathbb{C}$, we are able to find a suitable $D \subseteq \mathbb{C}^K$ such that $\mathbb{C}^K \cong \mathbb{E}^{K,\tr}(D)$. Again, one has to show appropriate forms of the conditions (i)---(iv) from Theorem~\ref{third}. The point of the Bays--Kirby method is we are able to hide the  transcendental number theory part of the Schanuel Property inside $D$, and so ignore it. In this way we are able to prove the quasiminimality of the powered fields $\mathbb{C}^K$ without characterising them all up to isomorphism.

\subsection{Outline of the paper}

The short section~\ref{terminology section} explains our terminology and notation for affine algebraic varieties, their linear counterparts, and the associated dimension notions.
The main technical objects of study, partial $K$-powered fields and their extensions, are introducted in section~\ref{definitions section}. 
Section~\ref{sec4} introduces the predimensions in the style of Hrushovski, which are the tool for expressing and using transcendence statements. 

The $K$-powers analogues of algebraicity and transcendence are explained in section~\ref{pregeometry section}, and we prove the first of the main technical steps in the quasiminimality proof, that purely powers-transcendental extensions can be amalgamated, using a lemma from stable group theory. 
In section~\ref{classextsec} we classify the finitely generated extensions of partial $K$-powered fields in terms of the locus of a good basis. This is possible due to a result of Zilber in Kummer theory. 

Section~\ref{a and e section} uses Hrushovski--Fra\"iss\'e amalgamation to build countable $K$-powered fields $\mathbb{F}^K(D_0)$ and $\mathbb{F}^{K, \tr}(D_0)$, which are then extended to continuum-sized quasiminimal $K$-powered fields $\mathbb{E}^K(D_0)$ and $\mathbb{E}^{K, \tr}(D_0)$ using the Shelah-Zilber excellence method.

In section~\ref{k-pow-closed fields sections}, we explain the $K$-powers closedness notion, and use a theorem of diophantine geometry known as ``weak Zilber--Pink'' to show that under a strong transcendence assumption, the analogue of Schanuel's conjecture, it implies an algebraic saturation property satisfied by $\mathbb{E}^K(D_0)$. 

In section~\ref{generic powers section} we put together the earlier work to characterise these  $\mathbb{E}^K(D_0)$ up to isomorphism by a short list of properties, which includes Theorem~\ref{third} as a special case. We then give several consequences, including Theorem~\ref{fourth}.

While sections~\ref{k-pow-closed fields sections} and~\ref{generic powers section} rely on a strong transcendence assumption which holds only for sufficiently generic subfields $K$ of $\mathbb{C}$, in section~\ref{quasiminimality section} we drop this assumption. The second main technical step towards the quasiminimality proof is to prove the algebraic saturation property for $\mathbb{E}^{K,\tr}(D_0)$ now only with tools from Ax's functional transcendence theorem in place of Schanuel's transcendence conjecture. We then complete the proof of Theorem~\ref{first}.

\subsection{Acknowledgements}
The second author would like to thank the Logic group at the University of Helsinki for their hospitality while this work was written up.

\section{Subspaces, loci, and dimension conventions}\label{terminology section}

There are two kinds of notions of \emph{subspace} related to a vector space which we make extensive use of, which model theorists characterize as semantic and syntactic notions, so to avoid ambiguity we briefly clarify the terminology and notation we will use. Let $K$ be a field, and $V$ a $K$-vector space.

By a \emph{($K$)-vector subspace} of $V$, as usual we mean a subgroup of $V$ closed under scalar multiplication.

Given a subset $A \subseteq V$, we write $\langle A \rangle_K$ for the \emph{$K$-span} of $A$ in $V$, the smallest $K$-vector subspace of $V$ containing $A$. 

Given subsets $A,B \subseteq V$, we write $\ldim_K(A/B)$, read as the \emph{$K$-linear dimension of $A$ over $B$}, to mean the cardinality of the smallest set $A_0$ such that $\langle A \cup B \rangle_K = \langle A_0 \cup B \rangle_K$; this is of course equal to the dimension of the quotient space $\gen{A \cup B}_K/\gen{B}_K$. As usual, we write $\ldim_K(A)$ for $\ldim_K(A/\emptyset)$.

Now fix $n \in \mathbb{N}$. By a \emph{$K$-linear subspace of $V^n$} we mean a $K$-vector subspace $L$ of $V^n$ given as the set of $n$-tuples $z$ from $V$ satisfying a matrix equation $Mz = 0$.
More generally, a \emph{$K$-affine subspace of $V^n$} is given as the solution set to a matrix equation $Mz = b$, for some matrix $M \in \textnormal{Mat}_{k \times n}(K)$ and some $b \in V^k$.
We say that a $K$-affine subspace $L$ is \emph{defined over a subset $B \subseteq V$} if it is defined by some equation $Mz=b$ with $b \in \langle B\rangle_K^k$.

The \emph{dimension} of a $K$-affine subspace $L$ of $V^n$ is defined to be $n-\rk(M)$, where $\rk(M)$ is the rank of the matrix $M$.

Now fix $a\in V^n$ and $B \subseteq V$. The \emph{$K$-affine locus of $a$ over $B$}, denoted by $K$-$\affloc(a/B)$, is the minimal $K$-affine subspace of $V^n$ which is defined over $B$ and contains $a$.
We can regard the finite tuple $a$ as a finite set, and then $\ldim_K(a/B)$ coincides with $\dim K$-$\affloc(a/B)$.

We have similar conventions for fields and algebraic varieties. 
Let $F$ be a field, and let $A, B$ be subsets of $F$. We write $\trdeg(A/B)$ for the cardinality of the smallest subset $A_0$ of $A$ such that every element of $F$ which is algebraic over $A\cup B$ is also algebraic over $A_0 \cup B$. Equivalently, in characteristic 0, $\trdeg(A/B)$ is the transcendence degree of the field extension $\mathbb Q(A\cup B) / \mathbb Q(B)$.

Given $a \in F^n$ and $B\subseteq F$, the \emph{algebraic locus} $\algloc(a/B)$ is the smallest Zariski-closed subset of $F^n$ which is defined over $B$ (or equivalently over the subfield of $F$ generated by $B$) and contains $a$.
We have $\dim \algloc(a/B) = \trdeg(a/B)$.

\section{$K$-powered fields}\label{definitions section}

In this section we introduce $K$-powered fields and the technical notions of a \textit{partial} or \textit{full} $K$-powered field, and discuss some basic facts about extensions.

\begin{defn}
Let $K$ be a field of characteristic 0. A \textit{$K$-powered field} consists of a $K$-vector space $V$, a field $F$ of characteristic 0, and a group homomorphism $$\exp:(V,+) \rightarrow (F^\times, \cdot).$$

We say the $K$-powered field is \textit{full} if $\exp$ is surjective and $F$ is algebraically closed.
\end{defn}

\begin{exa}\label{kpowexa}
For any subfield $K \subseteq \mathbb{C}$, the $K$-powered field $\mathbb{C}^K$ is the structure $$\mathbb{C}_{K\text{-VS}} \stackrel{\exp}{\longrightarrow} \mathbb{C}_{\textnormal{field}}$$  described in the Introduction. More generally, if $F_\exp$ is an exponential field and $K \subseteq F$ is a subfield, we construct a $K$-powered field $$F_{K\text{-VS}} \stackrel{\exp}{\longrightarrow} F_{\textnormal{field}}$$ where $F_{K\text{-VS}}$ is the reduct of $F_\exp$ to the $K$-vector space language, $F_\textnormal{field}$ is the reduct to the field language, and $\exp$ is the exponential map from $F_\exp$. 

The $K$-powered fields we will mostly be interested in in this paper have the form $\mathbb{B}^K$ and $\mathbb{C}^K$, where $\mathbb{B}_{\exp}$ is Zilber's exponential field and $\mathbb{C}_{\exp}$ is the complex exponential field.
\end{exa}

To classify extensions of $K$-powered fields it is useful to have a notion of \textit{partial} $K$-powered field, where the exponential is a partial function on $V$. 

\begin{defn}
A \textit{partial $K$-powered field} consists of a $K$-vector space $V$, a $\mathbb{Q}$-vector subspace $D$ of $V$, a field $F$ of characteristic 0, and a group homomorphism $\exp:D \rightarrow (F^\times)$ such that:
\begin{itemize}
    \item[1.] $\langle D \rangle_K=V$; and
    \item[2.] $\exp(D)$ generates $F$ as a field.
\end{itemize}

The \textit{kernel} $\ker_D$ of $D$ is the kernel of the exponential function. When it is clear what $D$ is, we will denote the kernel simply by $\ker$. A $K$-powered field $D$ has \textit{cyclic kernel} if $\ker_D$ is an infinite cyclic group.
\end{defn}

All the $K$-powered fields considered in this paper, such as $\mathbb{C}^K$ and $\mathbb{B}^K$, will have cyclic kernel.

A partial $K$-powered field is \textit{finitely generated} if the domain $D$ of the exponential is finite dimensional as a $\mathbb{Q}$-vector space.

\begin{exa}\label{standardbase}
Consider a 1-dimensional $\mathbb{Q}$-vector space $D$, generated by an element $\tau$. We take $V=D \otimes_\mathbb{Q} K$, $F=\mathbb{Q}(\sqrt{1})$ the field generated by all roots of unity, and define $\exp:D \rightarrow F$ by mapping $\frac \tau n$ to a primitive $n$-th root of unity, for each $n$, chosen so that for all $m_1$ and $m_2$ we have $$\exp \left( \frac \tau {m_1m_2} \right)^{m_2}=\exp \left( \frac{\tau }{m_1} \right).$$ This partial $K$-powered field is unique up to isomorphism. We will refer to this partial, 1-dimensional $K$-powered field $(D,V,\exp,F)$ as the \textit{standard base} and denote it by $SB^K$.
\end{exa}

It is possible to encode all the information of a partial $K$-powered field in the domain $D$ of the exponential map. 

\begin{defn}\label{language}
Let $K$ be a field of characteristic 0. Fix a $K$-vector space $V$ and an algebraically closed field $F$, also of characteristic 0. The language $\mathcal{L}^K$ is an expansion of the language of $\mathbb{Q}$-vector spaces by:

\begin{itemize}
\item[1.] A unary predicate $\ker$;
\item[2.] An $n$-ary predicate $L(z_1,\dots,z_n)$ for each $K$-linear subspace $L \leq V^n$;
\item[3.] An $n$-ary predicate $EW(w_1,\dots,w_n)$ for each algebraic variety $W \subseteq F^n$ defined and irreducible over $\mathbb{Q}$.
\end{itemize}
(The language does not depend on the choices of $V$ and $F$).
\end{defn}  

Given a partial $K$-powered field $(D,V,\exp,F)$ we may see the domain of exponentiation $D$ (which is a $\mathbb{Q}$-vector space by definition) as a structure in the language $\mathcal{L}^K$ as follows. The predicate $\ker$ is interpreted as the kernel of $\exp$; each predicate $L$ is interpreted as the intersection of the $K$-linear subspace $L \leq V^n$ with $D^n$; each predicate $EW$ is interpreted as the preimage under $\exp$ of the set of $F$-points of $W$.

Conversely, given the domain of exponentiation $D$ of a partial $K$-powered field seen as an $\mathcal{L}^K$-structure we can reconstruct the $K$-powered field.

\begin{nota}
   We will freely use $D$ (or $D_1,D_2,D'$...) to denote a partial $K$-powered field $(D,V,\exp,F)$. We will still denote by $F$ (or, according to the notation for the domain, $F_1,F_2,F'$...) the field generated by the quotient $D/\ker$, by $V$ (or $V_1,V_2,V'$...) the $K$-vector space that $D$ embeds in, by $\exp(z)$ the coset $z+\ker$ for some $z \in D$, and write ``$z \in L$'' or ``$\exp(z) \in W$'' rather than ``$L(z)$'' or ``$EW(z)$''. 
\end{nota}

\subsection{Extensions of partial $K$-powered fields}

\begin{defn}
Let $D_1$ be a partial $K$-powered field.

An \textit{extension} of $D_1$ is an $\mathcal{L}_K$-embedding $\phi$ of $D_1$ into a partial $K$-powered field $D_2$. Equivalently, it is an embedding of $(D_1,V_1,\exp_1,F_1)$ into $(D_2,V_2,\exp_2,F_2)$ consisting of a $K$-linear embedding $\phi_V:V_1 \hookrightarrow V_2$ and a field embedding $\phi_F:F_1 \hookrightarrow F_2$ such that $\phi_V(D_1) \subseteq D_2$ and $\phi_F \circ \exp_1=\exp_2 \circ \phi_V$.

If $\phi$ is an inclusion, we say that $D_1$ is a \textit{$K$-powered subfield} of $D_2$. We denote this by $D_1 \leq D_2$. If $\ker_{D_2}=\ker_{D_1}$, we say that the extension \textit{preserves the kernel}, or that it is \textit{kernel-preserving}.
\end{defn}

\begin{exa}
    Observe that a partial $K$-powered field has cyclic kernel if and only if it is a kernel-preserving extension of the standard base $SB^K$ of Example \ref{standardbase}.
\end{exa}

An extension of partial $K$-powered fields $D_1 \leq D_2$ is \textit{finitely generated} if $D_2$ is finitely generated over $D_1$ as a $\mathbb{Q}$-vector space. A \textit{basis} for the extension $D_2$ of $D_1$ is a $\mathbb{Q}$-linear basis for $D_2$ over $D_1$.

We can use bases to determine the isomorphism type of an extension of $K$-powered fields. 

\begin{defn}
    Let $D$ be a partial $K$-powered field, $z \in D^n$, $A \subseteq D$. The $K$-\textit{powers locus} of $z$ over $A$ is the pair $(K$-$\affloc(z/A), \algloc(\exp(z)/\exp(A)))$. We denote it by $\Loc(z/A)$.
\end{defn}

The following is straightforward.

\begin{lem}\label{fingeniso}
Let $D_1$ and $D_2$ be finitely generated extensions of a partial $K$-powered field $D$, and let $b_1 \in D_1^n$ be a basis for $D_1$ over $D$.

Then $D_1$ and $D_2$ are isomorphic over $D$ if and only if there is a basis $b_2 \in D_2^n$ such that for each $m \in \mathbb{N} \setminus \{0\}$, we have $\Loc \left(\frac{b_1}{m}/D\right)=\Loc \left( \frac{b_2}{m}/D \right)$ . \hfill $\qedsymbol{}$
\end{lem}

In Section \ref{classextsec} we will see that in cases of interest it is actually sufficient to consider the loci of $\frac{b}{m}$ up to some finite $m$, and therefore by replacing $b$ by $\frac{b}{m!}$ we can consider only $m=1$.

\section{Predimensions}\label{sec4}

In this section we introduce the predimension $\delta^K$ on a $K$-powered field and we use it to define \textit{strong extensions} of $K$-powered fields. We show that these classes of extensions satisfy an amalgamation property. 

We then compare the predimension on $\mathbb{B}^K$ with the exponential predimension on Zilber's exponential field $\mathbb{B}_\exp$, and use this comparison to show that if $K$ is a subfield of $\mathbb{B}_\exp$ of finite transcendence degree then $\mathbb{B}^K$ has a finitely generated strong substructure. 
 
\subsection{The predimension on a $K$-powered field}

\begin{defn}
Let $D$ be a partial $K$-powered field, $A,B$ subsets of $D$ with $\ldim_\mathbb{Q}(A/B)$ finite.

We define the \textit{predimension} of $A$ over $B$ as $$\delta^K(A/B)=\ldim_K(A/B)+\trdeg(\exp(A)/\exp(B))-\ldim_\mathbb{Q}(A/B).$$

We define $\delta^K(A):=\delta^K(A/\varnothing)$.
\end{defn}

When the field of powers is clear we just write $\delta$ instead of $\delta^K$. We will frequently consider the predimension of a tuple $z \in D^n$ rather than of a set; this is defined in the obvious way. 

Note that if $D_1 \leq D_2$ is a finitely generated extension of partial $K$-powered fields, then $\delta(D_2/D_1)$ is equal to $\delta(z/D_1)$ where $z$ is a basis for $D_2$ over $D_1$.

\begin{lem}
The predimension function $\delta^K$ satisfies the following properties:

\begin{itemize}
\item[1.] Finite character: if $z \in D^n$ and $A \subseteq D$, there is a finite subset $A_0 \subseteq A$ such that $\delta(z/A)=\delta(z/A_0)$.
\item[2.] Addition formula: given $X \subseteq Y \subseteq Z \subseteq D$, $$\delta(Z/X)=\delta(Z/Y)+\delta(Y/X).$$
\item[3.] Submodularity: given $A_1,A_2 \subseteq D$, $$\delta(A_1 \cup A_2/A_1) \leqslant \delta (A_2 / \langle A_1 \rangle_\mathbb{Q} \cap \langle A_2 \rangle_\mathbb{Q}).$$ 
\end{itemize}
\end{lem}

\begin{proof}
Straightforward, by the same argument as \autocite[Lemma 4.2]{BK18}.
\end{proof}

\subsection{Strong extensions}

As usual in Hrushovski--Fra\"iss\'e-type constructions, the predimension function is used to define a notion of \textit{strong extension}.

\begin{defn}
An extension $D_1 \leq D_2$ of partial $K$-powered fields is \textit{strong} if it is kernel-preserving and for all finite tuples $z$ in  $D_2$, we have $\delta(z/D_1) \geqslant 0$. We denote this by $D_1 \treq D_2$ and say $D_1$ is a \textit{strong $K$-powered subfield} of $D_2$.
\end{defn}

It is straightforward to verify that the composite of strong extensions is strong.

\begin{defn}[Free amalgam]
Let $D_2$ and $D_3$ be partial $K$-powered fields, extending the full $K$-powered field $D_1$.

Let $V_4$ be the $K$-vector space $V_2 \oplus_{V_1} V_3$, and $D_4$ the $\mathbb Q$-vector subspace of $V_4$ generated by $D_2$ and $D_3$. Since $F_1$ is algebraically closed, there is up to isomorphism a unique amalgam $F_4$ of $F_2$ and $F_3$ such that $F_2$ is linearly disjoint from $F_3$ over $F_1$. Extend $\exp$ to $D_4$ by setting $\exp(z_2+z_3)=\exp(z_2) \cdot \exp(z_3)$.

We call the partial $K$-powered field $D_4$ the \textit{free amalgam} of $D_2$ and $D_3$ over $D_1$.
\end{defn}

\begin{prop}\label{fullamalgam}
Let $D_2$ and $D_3$ be partial $K$-powered fields, extending the full $K$-powered field $D_1$. Let $D_4$ be the free amalgam of $D_2$ and $D_3$ over $D_1$, and suppose $D_1 \treq D_2$. Then $D_3 \treq D_4$. If also $D_1 \treq D_3$, then $D_1 \treq D_4$.
\end{prop}

\begin{proof}
First we note that $\ker_{D_4}=\ker_{D_3}$: if $\exp(z_2+z_3)=1$ then $\exp(z_2) = \frac{1}{\exp(z_3)} \in F_3$, so $\exp(z_2) \in F_1$. Since $D_1 \treq D_2$ and $D_1$ is full, this implies that $z_2 \in D_1$. Hence $z_2+z_3 \in D_3$.

Let now $z_4$ be a finite tuple in $D_4$; then $z_4=z_2+z_3$ for some finite tuples $z_2$ and $z_3$. Then $\delta(z_4/D_3)=\delta(z_2/D_3)$. Since the amalgam is free, we have $\delta(z_2/D_3)=\delta(z_2/D_1)$, which is non-negative because $D_1 \treq D_2$. So indeed $D_3 \treq D_4$.

The last statement follows from the composite of strong extensions being strong.
\end{proof}

\subsection{Hulls}

\begin{lem}\label{strongint}
Suppose $D$ is a partial $K$-powered field, $J$ is a set, and for each $j \in J$, $D_j$ is a strong $K$-powered subfield of $D$.

Then $D_J:=\bigcap_{j \in J} D_j$ is also strong in $D$.
\end{lem}

\begin{proof}
By the same argument as \autocite[Lemma 4.5]{BK18}.
\end{proof}

\begin{defn}
Let $D$ be a partial $K$-powered field, $A \subseteq D$. The \textit{hull} of $A$ in $D$ is $$\hull{A}^D_K=\bigcap \{X \treq D \mid A \subseteq X\}.$$
\end{defn}

Lemma \ref{strongint} shows that the hull of a set in $D$ is a strong $K$-powered subfield of $D$. When the field of exponents and the ambient space in which we take the hull are clear, we drop this from the notation and write $\hull{A}$ rather than $\hull{A}^D_K$.

\begin{lem}
The hull operator has finite character: if $D$ is a $K$-powered field and $X \subseteq D$, then $$\hull{X}=\bigcup_{X_0 \subseteq X \textnormal{ finite}} \hull{X_0}.$$
\end{lem}

\begin{proof}
By the same argument as \autocite[Lemma 4.7]{BK18}.
\end{proof}

\subsection{Strong substructures of $\mathbb{B}^K$}

In this subsection we study strong substructures of the $K$-powered field $\mathbb{B}^K$, coming from Zilber's exponential field as explained in Example \ref{kpowexa}. To do so, we compare the predimension $\delta^K$ with the notion of predimension used in exponential fields. We will explain the properties of $\mathbb{B}_\exp$ we need as we use them.

\begin{defn}
Let $F$ be an exponential field. The \textit{exponential predimension} $\delta^\exp$ on $F$ is the predimension function defined, for subsets $A,B \subseteq F$ with $\ldim_{\mathbb{Q}}(A/B)$ finite, by $$\delta^\exp(A/B)=\trdeg(A,\exp(A)/B,\exp(B))-\ldim_\mathbb{Q}(A/B).$$

As in the case of the predimension $\delta^K$, we will frequently refer to the predimension of a tuple $z$ in the exponential field $F$. 

For a set $B \subseteq F$, we write $B \treq_{\exp} F$ to denote the fact that $\delta^{\exp}(z/B) \geqslant 0$ for all finite tuples $z$ in $F$, and say that $B$ is \textit{exponentially strong} in $F$.

The \textit{exponential hull} $\hull{B}^F_\exp$ of a set $B$ is the smallest $\mathbb{Q}$-vector subspace of $F$ containing $B$ that is exponentially strong in $F$ (and it always exists).
\end{defn}

In this case too we will omit the field $F$ from the notation if it is clear from the context. 

\begin{lem}\label{hullinclusion}
Let $F$ be an exponential field, $K \subseteq F$ a subfield, $X \subseteq F$. Then $\hull{X}_K \subseteq \hull{X,K}_\exp$.
\end{lem}

\begin{proof}
Let $D=\lceil X,K \rceil_{\exp}$: we have to prove that $D$ is strong in $F^K$ in the sense of $K$-powered fields. So let $z$ be a tuple in $F$. Then
\begin{align*}
\delta^K(z/D) &= \ldim_K(z/D) + \trdeg(\exp(z)/\exp(D)) -\ldim_\mathbb{Q}(z/D) \\
& \geqslant \trdeg(z/D)+\trdeg(\exp(z)/\exp(D)) - \ldim_\mathbb{Q}(z/D) \\
& \geqslant \trdeg(z,\exp(z)/D,\exp(D))-\ldim_\mathbb{Q}(z/D) \\ 
& = \delta^{\exp}(z/D) \\
& \geqslant 0
\end{align*}
as required.
\end{proof}

We will use this with $F=\mathbb{B}_\exp$ and $X=\{\tau\}$, where $\tau$ is a generator of the kernel in $\mathbb B_\exp$. Strong extensions are by definition kernel-preserving, and therefore $\tau$ is contained in any strong $K$-powered subfield of $\mathbb B^K$. Hence the hull of the empty set, that is, the minimal strong substructure of $\mathbb{B}_\exp$, coincides with the hull of $\tau$. 

Recall that by construction (see \autocite[Section 5]{Zil05}, \autocite[Section 9]{BK18}) $\mathbb B_\exp$ satisfies the Schanuel property, that is, for every finite $A \subseteq \mathbb B_\exp$ we have $\delta^\exp(A) \geqslant 0$.

\begin{lem}\label{fingenstrong}
Let $K \subseteq \mathbb{B}_\exp$ be a subfield of transcendence degree $d$ for some $d \in \mathbb{N}$. Then $\mathbb{B}^K$ has a finitely generated strong partial $K$-powered subfield (in particular $\lceil \tau \rceil_K$ is finitely generated).
\end{lem}

\begin{proof}
Let $z$ be a finite tuple in $\mathbb{B}_\exp$. Then we have:
\begin{align*}
\delta^K(z) &=\ldim_K(z)+\trdeg(\exp(z))-\ldim_{\mathbb{Q}}(z) \\
& \geqslant \trdeg(z/K) +\trdeg(\exp(z)) - \ldim_{\mathbb{Q}}(z)\\
& \geqslant \trdeg(z) - d + \trdeg(\exp(z)) - \ldim_{\mathbb{Q}}(z) \\
& \geqslant \trdeg(z,\exp(z)) - \ldim_{\mathbb{Q}}(z) -d \\
& = \delta^{\exp}(z) - d \\
& \geqslant -d.
\end{align*}

Thus there is a tuple $a$ containing $\tau$ such that $\delta^K(a)$ is minimal. Then for any $z$ we have $\delta^K(z/a) \geqslant 0$, so $\langle a \rangle_\mathbb{Q}$ is a finitely generated strong $K$-powered subfield.
\end{proof}

\begin{rem}\label{smallhull}
    By \autocite[Theorem 1.3]{BK18}, there are many choices of $K$ for which $\hull{\tau}_K \cong SB^K$; for instance, this holds when $K=\mathbb{Q}(\lambda)$ for all but countably many $\lambda$'s in $\mathbb{B}_\exp$.
\end{rem}

\section{The $K$-powers-closure pregeometry}\label{pregeometry section}

As usual in Hrushovski-Fra\"iss\'e amalgamation constructions, we use the predimension to define a pregeometry on $K$-powered fields. We introduce two kinds of extensions of $K$-powered fields: the \textit{powers-algebraic} and \textit{powers-transcendental} extensions, analogous to algebraic and purely transcendental extensions of fields. Finally, we show that on $\mathbb{C}^K$ the pregeometry satisfies the countable closure property using the corresponding property of exponential-algebraic closure.

\subsection{The pregeometry on a $K$-powered field}

\begin{defn}
Let $D$ be a partial $K$-powered field, and let $A \subseteq D$. The \textit{$K$-powers-closure of $A$ in $D$}, denoted $\pcl_K^D(A)$, is the smallest $K$-powered subfield of $D$ which contains $A$ such that if $B \subseteq D$ and $\delta(B/\pcl_K^D(A)) \leqslant 0$ then $B \subseteq \pcl_K^D(A)$.
\end{defn}

If the field of exponents is clear then we only write $\pcl^D(A)$.

\begin{prop}
Let $D$ be a partial $K$-powered field. Then $$\pcl_K^D(A)=\bigcup \{z \in D \mid \exists n \in \mathbb{N}, z' \in D^n : \delta(z,z'/\hull{A}_K^D)=0 \}$$ and $\pcl^D$ is a pregeometry on $D$. Moreover, if $D$ is full and $A \subseteq D$, then the powers-closure $\pcl(A)$ is a full $K$-powered subfield of $D$.
\end{prop}

\begin{proof}
The description of $\pcl_K^D(A)$ is straightforward. Finite character is by the same argument as \autocite[Lemma 4.12]{BK18}, and exchange is proved easily. The ``moreover'' statement is by the same argument as \autocite[Lemma 4.10(2)]{BK18}.
\end{proof}

As with all pregeometries, we can associate to powers-algebraic closure a dimension function.

\begin{defn}
    We call the dimension function associated to $\pcl^D$ \textit{powers-transcendence degree} and denote it by $\pdim^D$.

    If $z \in D$ and $A \subseteq D$ we say that $z$ is \textit{powers-algebraic over $A$ in $D$} if $\pdim^D(z/A)=0$, and \textit{powers-transcendental over $A$ in $D$} otherwise. 
\end{defn}

We have defined $\pcl_K^D$ as a closure operator on the covering sort $D$ of a partial $K$-powered field. For a full $K$-powered field, we extend it to a closure operator on the disjoint union of the covering and field sorts $D \sqcup F$, by also closing under $\exp$ and $\exp^{-1}$. If we now restrict to the field sort $F$ we get a pregeometry there which we denote by $\pcl_K^F$.

In the case $K=\mathbb{Q}$, we have $\delta^{\mathbb{Q}}(z/A)=\trdeg(\exp(z)/\exp(A))$ and it is easy to see that $\pcl_\mathbb{Q}^F$ is equal to algebraic closure on $F$, that is, given $b \in F$ and a subfield $A \subseteq F$, we have $b \in \pcl_{\mathbb{Q}}^F(A)$ if and only if $b$ is a zero of a non-trivial polynomial $p(X) \in A[X]$.

Now consider $K \neq \mathbb{Q}$ and suppose $b \in F$ is a zero of a ``$K$-powered polynomial'' over a subfield $A \subseteq F$, that is, there is a non-zero polynomial $p(X_1,\dots,X_r) \in A[X_1,\dots,X_r]$, there are $\mathbb{Q}$-linearly independent $\lambda_1,\dots,\lambda_r \in K$, and there are $a_1,\dots,a_r \in \exp^{-1}(b)$ such that $p(\exp(\lambda_1a_1),\dots,\exp(\lambda_ra_r))=0$. We also write this informally as $p(b^{\lambda_1},\dots,b^{\lambda_r})=0$.

Then one can check easily that $\delta^K(\lambda_1a_1,\dots,\lambda_ra_r/\hull{\exp^{-1}(A)}_K^D)=0$ and so $b \in \pcl_K^F(A)$.

However, we can also have $b \in \pcl_K^F(A)$ without satisfying a $K$-powers polynomial over $A$, due to a lack of elimination theory for $K$-powered polynomial.

For example, if we take $\lambda_1,\lambda_2,c_1,\dots,c_8$ to be a sufficiently generic 10-tuple in $\mathbb{B}_\exp$, for example exponentially-algebraically independent, then we can find $a_1,a_2 \in \mathbb{B}_\exp$ such that $$c_1\exp(a_1)+c_2\exp(\lambda_1a_1)+c_3\exp(a_2)+c_4\exp(\lambda_2a_2)=c_5\exp(a_1)+c_6\exp(\lambda_1a_1)+c_7\exp(a_2)+c_8\exp(\lambda_2a_2)=0.$$

So taking $K=\mathbb{Q}(\lambda_1,\lambda_2)$ and $C=\mathbb{Q}(c_1,\dots,c_8)$ we see that $a_1,a_2$ satisfy two different $K$-powered polynomials in two variables, and one can see that $a_1,a_2 \in \pcl_K^{\mathbb{B}_\exp}(C)$.

However, using the arguments of \autocite[Proposition 7.5]{Kir13}, one can check that neither $a_1$ nor $a_2$ satisfies any one-variable $K$-powered polynomial over $C$.

The predimension and the pregeometry are related in the usual way.

\begin{lem}\label{predimpregeo}
Let $D$ be a partial $K$-powered field. Suppose $D_1 \treq D$ and $D_2$ is a finitely generated extension of $D_1$ contained in $D$. Then:

\begin{itemize}
    \item[$(1)$] $D_2 \treq D$ if and only if $\pdim^D(D_2/D_1)=\delta(D_2/D_1)$;
    \item[$(2)$] $\pdim^D(D_2/D_1)=\min \{ \delta(D'/D_1) \mid D_2 \subseteq D' \subseteq D \}$; 
    \item[$(3)$] $\hull{D_2}^D$ is finitely generated over $D_2$.
\end{itemize}
\end{lem}

\begin{proof}
These are standard facts about the pregeometry attached to a predimension. A proof of the first two assertions follows the same argument as \autocite[Lemma 4.15]{BK18} (although in this setting there is no need to apply \autocite[Lemma 4.13]{BK18}.)

The last statement is a consequence of the first two: the extension finitely generated over $D_2$ by the elements of a $\pcl$-basis of $\pcl^D(D_2)$ over $D_1$ is strong in $D$.
\end{proof}

Unlike the usual transcendence and algebraicity, the ambient $D$ does matter in general. However, the following lemma, standard for predimension constructions, shows that as long as extensions are strong we can forget the ambient $D$.

\begin{lem}\label{strong ext preserve ktd}
   Let $D_0 \treq D_1 \treq D_2$ be partial $K$-powered fields. For every tuple $a$ in $D_1$, $\pdim^{D_1}(a/D_0)=\pdim^{D_2}(a/D_0)$. 
\end{lem}

\begin{proof}
    Follows from Lemma \ref{predimpregeo}(2).
\end{proof}

\subsection{Powers-algebraic extensions}

\begin{defn}\label{powalgdef}
A strong extension $D_1 \treq D_2$ of partial $K$-powered fields is \textit{powers-algebraic} if every element in $D_2$ is powers-algebraic over $D_1$; equivalently, if for all finite tuples $z$ in $D_2$ there is a finite tuple $z'$ in $D_2$ extending $z$ such that $\delta(z'/D_1)=0$.
\end{defn}

\begin{prop}
    Let $D_1$ be a full $K$-powered field and let $D_2$ and $D_3$ be powers-algebraic extensions of $D_1$. The free amalgam $D_4$ of $D_2$ and $D_3$ over $D_1$ is a powers-algebraic extension of $D_1$.
\end{prop}

\begin{proof}
    Let $z_2+z_3$ be a finite tuple in $D_4$, and let $D_2'=D_2  \cap \langle D_1, z_2+z_3 \rangle_{\mathbb{Q}}$.
    \begin{align*}
        \delta(z_2+z_3/D_1) &= \delta (z_2+z_3/D_2')+\delta(D_2' /D_1)    \\
                            & = \delta(z_3/D_2')  \\
                            & = \delta(z_3/D_1)     \,\,\,\,\,\,\,\,\,\,\,\,\,\,\,\,\,\,\,\,\,\,\,\,\,\,\,\,\,\,\,\,\,\,\,\,\,\,\,\,\,\,\,\,\,\,\,\,\,\,\,\,\,\,\,\textnormal{(by freeness of the amalgam)} \\
                            &=0. \qedhere
                            \end{align*}
                            \end{proof}

\subsection{Purely powers-transcendental extensions}

\begin{defn}
Let $D_1$ be a $K$-powered subfield of $D_2$. We say $D_1$ is \textit{$K$-powers-closed} in $D_2$ if for all finite tuples $z$ in $D_2$ we have that either $\delta(z/D_1)>0$, or $z$ is in $D_1$. Equivalently, if $D_1=\pcl^{D_2}(D_1)$.

In this case we say that $D_2$ is a \textit{purely powers-transcendental extension} of $D_1$. We denote this by $D_1 \treqcl D_2$.
\end{defn}

We next prove an amalgamation result for purely powers-transcendental extensions of a fixed, full $K$-powered field. This is one of the two key technical steps in the Bays-Kirby method which allows us to prove quasiminimality of $\mathbb{C}^K$ without Schanuel's Conjecture. 

\begin{prop}\label{pptamal}
Let $D_0$ be a full, countable $K$-powered field, let $D_1, D_2,$ $D_3$ be full, purely powers-transcendental extensions of $D_0$ with $D_1 \treq D_2$ and $D_1 \treq D_3$, and let $D_4$ be the free amalgam of $D_2$ and $D_3$ over $D_1$. Then $D_4$ is a purely powers-transcendental extension of $D_0$.
\end{prop}

We are going to use two versions of a technical result on stable groups due to Ziegler.

\begin{lem}[{\autocite[Theorem 1]{Zie06}}]\label{stablegroups}
\textbf{Vector space version.}
Let $V$ be a $K$-vector space, $V' \leq V$ a vector subspace, and let $v_1,v_2,v_3 \in V^n$. Suppose that $v_1+v_2+v_3=0_{V^n}$ and that $$\ldim_K(v_i/V'v_j)=\ldim_K(v_i/V')$$ for each $i,j$ in $\{1,2,3\}$ with $i \neq j$.

Then there is a $K$-linear subspace $L \leq V^n$ such that $v_1,v_2,v_3$ are generic points of $V'$-cosets of $L$; in particular, $\ldim_K(v_i/V')=\dim L$ for $i=1,2,3$.

\textbf{Algebraic group version.}
Let $H$ be a commutative algebraic group defined over an algebraically closed field $F$, and let $h_1,h_2,h_3 \in H$. Suppose that $h_1+h_2+h_3=0_H$, and that $$\trdeg(h_i/Fh_j)=\trdeg(h_i/F)$$ for each $i,j$ in $\{1,2,3\}$ with $i \neq j$.

Then there is an algebraic subgroup $G$ of $H$ such that $h_1,h_2,h_3$ are generic points of $F$-cosets of $G$; in particular, $\trdeg(h_i/F)=\dim G$ for $i=1,2,3$.
\end{lem}

\begin{proof}[Proof of Proposition \ref{pptamal}]
We will show that $\delta(b/D_0)>0$ for every tuple $b$ in $D_4$. Assume first that $\langle D_0, b \rangle_\mathbb{Q} \cap D_2 \neq D_0$. Then we have 
\begin{align*}
 \delta(b / D_0) & =\delta(b/\langle D_0,b \rangle_\mathbb{Q} \cap D_2)+ \delta(\langle D_0,b \rangle_\mathbb{Q} \cap D_2/D_0) \geqslant \\
 &\geqslant \delta(b/D_2) + \delta(\langle D_0,b \rangle_\mathbb{Q} \cap D_2/D_0) \geqslant \,\,\,\,\,\,\,\,\,\,\,\,\,\,\,\,\,\, \textnormal{(by submodularity of }\delta\textnormal{)} \\  &\geqslant 0+1 > 0
\end{align*}
because $D_2 \treq D_4$ and $\langle D_0,b \rangle_\mathbb{Q} \cap D_2$ is contained in $D_2$ which is a purely powers-transcendental extension of $D_0$. 

Now we assume $\langle D_0, b \rangle_\mathbb{Q} \cap D_2=D_0$ and, symmetrically, $\langle D_0, b \rangle_\mathbb{Q} \cap D_3=D_0$.

By definition of free amalgam, there are $b_2 \in D_2^n$ and $b_3 \in D_3^n$ such that $b=b_2+b_3$. We denote by $D_2'$ and $D_3'$ the partial $K$-powered fields $\langle D_0,b_2 \rangle_\mathbb{Q}$ and $\langle D_0, b_3 \rangle_\mathbb{Q}$ respectively. We claim that $D_2 \cap D_3'=D_2' \cap D_3=D_0$.

To see this, let $v \in D_2 \cap D_3'$. Since it is in $D_3'$, we may write it as $v_0 + \sum_{i=1}^n q_ib_3^i$ for some $v_0 \in D_0$, $q_1,\dots,q_n \in \mathbb{Q}$, where $b_3=(b_3^1,\dots,b_3^n)$. Let then $u_3=v-v_0=\sum_{i=1}^n q_i b_3^i$, and define analogously $u_2=\sum_{i=1}^n q_i b_2^i$ and $u=u_2+u_3=\sum_{i=1}^n q_i b^i$. Since $b_2 \in D_2^n$, we have $u=u_2+u_3 \in D_2$: but we have assumed that $\langle D_0, b \rangle_\mathbb{Q} \cap D_2=D_0$, and therefore $u \in D_0$. Therefore all the $q_i$'s are actually 0, and we have $v=v_0 \in D_0$. Hence $D_2 \cap D_3'=D_0$; the same argument proves that $D_2' \cap D_3=D_0$. 

Consider now the $K$-powered field $C=\langle D_2', b \rangle_\mathbb{Q}$ and note that since $b=b_2+b_3$, $C=\langle D_2', b_3 \rangle_\mathbb{Q}=D_2'+D_3'$. Consider the following extensions of $\mathbb{Q}$-vector spaces:

$$
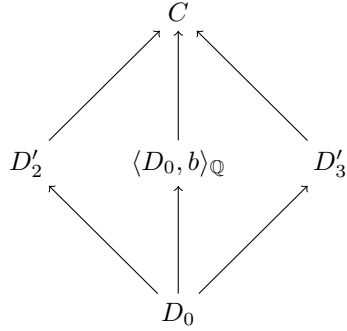
\begin{tikzpicture} [node distance=2cm]
\node (a) at (0,0) {$D_0$};
\node (b) at (-2,2) {$D_2'$};
\node (c) at (0,2) {$\langle D_0, b \rangle_\mathbb{Q}$};
\node (d) at (2,2) {$D_3'$};
\node (e) at (0,4) {$C$};

\draw [->] (a) -- (b);
\draw [->] (a) -- (c);
\draw [->] (a) -- (d);
\draw [->] (b) -- (e);
\draw [->] (c) -- (e);
\draw [->] (d) -- (e);
\end{tikzpicture}$$

By assumption, $\ldim_\mathbb{Q}(b/D_0)=n$. By modularity of $\mathbb{Q}$-linear dimension applied to the half-square diagrams, we have $$\ldim_\mathbb{Q}(C/D_3')=\ldim_\mathbb{Q}(C/D_2') = \ldim_\mathbb{Q}(b/D_0)=n$$ and $$\ldim_\mathbb{Q}(D_2'/D_0)=\ldim_\mathbb{Q}(C/\langle D_0, b \rangle_\mathbb{Q})=\ldim_\mathbb{Q}(D_3'/D_0)$$ while applying modularity to the full square we get
\begin{align*}
    \ldim_\mathbb{Q}(D_2'/D_0) &= \ldim_\mathbb{Q}(C/D_3') \\
    \ldim_\mathbb{Q}(D_3'/D_0) &= \ldim_\mathbb{Q}(C/D_2'). 
\end{align*}
Combined, these imply that $b_2$ and $b_3$ are both $\mathbb{Q}$-linearly independent over $D_0$, and in fact over $D_1$ since $D_1 \cap D_2' \subseteq D_3 \cap D_2'=D_0$ (and similarly for $D_1 \cap D_3'$).

We now consider some inequalities.
\begin{align}
    \ldim_K(b/D_1) \geqslant \ldim_K(b/D_1b_3) = \ldim_K(b_2/D_1b_3) =\ldim_K(b_2/D_1) \label{one}
\end{align} where the first equality holds because $b=b_2+b_3$ and the second one because $D_2$ and $D_3$ are amalgamated freely in $D_4$, and therefore $b_2$ and $b_3$ are independent over $D_1$ in the sense of $K$-vector spaces. With similar reasoning, we obtain 
\begin{align}
    \ldim_K(b/D_1) & \geqslant \ldim_K(b_3/D_1) \label{two} \\
    \trdeg(\exp(b)/F_1) & \geqslant \trdeg(\exp(b_2)/F_1) \label{three} \\
    \trdeg(\exp(b)/F_1) & \geqslant \trdeg(\exp(b_3)/F_1). \label{four}
\end{align}

Moreover, we have 
\begin{align}
    \ldim_K(b/D_0) &\geqslant \ldim_K(b/D_1) \label{five}\\ 
    \trdeg(\exp(b)/F_0) &\geqslant \trdeg(\exp(b)/F_1) \label{six}
\end{align} simply because linear dimension and transcendence degree do not increase when we extend the base.

Now if one of (\ref{five}) and (\ref{six}) is a strict inequality, we obtain 
\begin{align*}
    \delta(b/D_0) &=\ldim_K(b/D_0)+\trdeg(\exp(b)/F_0) - \ldim_\mathbb{Q}(b/D_0) \\
    &> \ldim_K(b/D_1)+\trdeg(\exp(b)/F_1) - n \\
    & \geqslant \ldim_K(b_2/D_1)+\trdeg(\exp(b_2)/F_1) -\ldim_\mathbb{Q}(b_2/D_1) \\
    & =\delta(b_2/D_1) \geqslant 0
\end{align*} and therefore $\delta(b/D_0)>0$, as we wanted. 

Assume then that (\ref{five}) and (\ref{six}) are both equalities, and that one of (\ref{one}) and (\ref{three}) is strict. Again, we obtain
\begin{align*}
    \delta(b/D_0) &= \ldim_K(b/D_0)+\trdeg(\exp(b)/F_0) - \ldim_\mathbb{Q}(b/D_0) \\
    & = \ldim_K(b/D_1) + \trdeg(\exp(b)/F_1) -n \\
    & > \ldim_K(b_2/D_1) + \trdeg(\exp(b_2)/F_1) -n \\
    & = \delta(b_2/D_1) \geqslant 0
\end{align*} and therefore again $\delta(b/D_0)>0$. If (\ref{five}) and (\ref{six}) are equalities, and one of (\ref{two}) and (\ref{four}) is strict, we obtain the same replacing $b_2$ with $b_3$. Thus we assume that none of the six inequalities is strict.

In this case, we obtain that
\begin{align*}
    \ldim_K(b/D_1)&=\ldim_K(b/D_1b_2) \\
    \ldim_K(b/D_1)&=\ldim_K(b/D_1b_3) \\
    \trdeg(\exp(b)/F_1)&=\trdeg(\exp(b)/D_1\exp(b_2)) \\
    \trdeg(\exp(b)/F_1)&=\trdeg(\exp(b)/D_1\exp(b_3))
\end{align*} which together with the fact that $D_4$ is a free amalgam imply that $b,b_2,b_3$, and $\exp(b),\exp(b_2),\exp(b_3)$ satisfy the assumptions of the two versions of Lemma \ref{stablegroups}. Thus there are a $K$-linear subspace $L \leq V_4^n$ and an algebraic subgroup $H \subseteq ((F_4^\alg)^\times)^n$ such that $b,b_2,b_3$ are generic points of $D_1$-translates of $L$ and $\exp(b),\exp(b_2),\exp(b_3)$ are generic points of $F_1$-cosets of $H$. Thus we have $\dim L=\ldim_K(b/D_1)$ and $\dim H=\trdeg(\exp(b)/F_1)$, and as we are assuming that (\ref{five}) and (\ref{six}) are equalities this says that $\dim L = \ldim_K(b/D_0)$ and $\dim H=\trdeg(\exp(b)/F_0)$; in other words, $$\delta(b/D_0)=\dim L + \dim H -n.$$

Since $\ldim_\mathbb{Q}(b/D_0)=n$, we must have $b \notin D_0^n$, so $\ldim_K(b/D_0)>0$. 

Assume now that $m_1,\dots,m_n \in \mathbb{Z}$ satisfy $$\exp(b^1)^{m_1} \cdots \exp(b^n)^{m_n} =c$$ for some $c \in F_0$; then $c=\exp(a)$ for some $a \in D_0$ (since $D_0$ is full) and therefore
\begin{align*}
    \exp(m_1b^1 + \dots +m_nb^n)&=\exp(a) \\
    \exp(m_1b^1 + \dots +m_nb^n-a)&=1
\end{align*} so $m_1b^1 + \dots +m_nb^n-a \in \ker_{D_2+D_3}$. All extensions are strong, and therefore kernel-preserving, so $\ker_{D_2+D_3}=\ker_{D_0} \subseteq D_0$: hence, this implies that $m_1b^1 + \dots +m_nb^n \in D_0$, which can only hold for $m_1=\dots=m_n=0$, again because $\ldim_\mathbb{Q}(b/D_0)=n$. Therefore, $\exp(b)$ is not contained in any $F_0$-coset of a proper algebraic subgroup of $((F_4^\alg)^\times)^n$; hence $\dim H=n$. So we have
\begin{align*}
    \delta(b/D_0) &= \dim L + \dim H -n\\
    &> 0 +n-n=0.
\end{align*} This completes the proof.
\end{proof}

\subsection{The countable closure property}\label{subsecpreoncb}

\begin{defn}
Let $A$ be a set with a pregeometry $\cl$. We say $(A,\cl)$ has the \textit{countable closure property} if for all finite $X \subseteq A$, $\cl(X)$ is countable.
\end{defn}

We sum up in the following statement the main results obtained by the second author in \cite{Kir10}.

\begin{thm}[{\autocite[Theorems 1.1 and 1.2]{Kir10}}]\label{ecl}
Any exponential field $F$ has a pregeometry operator $\ecl^F$, with associated dimension $\edim^F$, such that every $\ecl^F$-closed subset of $F$ is an exponential subfield; if $C$ is such a set, then every $z \in F^n$ satisfies $$\delta^\exp(z/C)=\trdeg(z,\exp(z)/C)-\ldim_\mathbb{Q}(z/C) \geqslant \edim^F(z/C).$$
\end{thm}

The fact that $\mathbb{B}_\exp$ has the countable closure property is true by construction of $\mathbb{B}_{\exp}$ (see for example \autocite[Lemma 5.11]{Zil05b}). The countable closure property for $\mathbb{C}_\exp$ is \autocite[Lemma 5.12]{Zil05b}.

Consider now a countable subfield $K$ of an exponential field $F$ with the countable closure property. For any exponential subfield $C \subseteq F$ containing $K$ we may consider the $K$-powered field $C^K$.

\begin{prop}\label{eclfull}
Let $F$ be an exponential field which is algebraically closed, such that $\exp$ is surjective, with the countable closure property for $\ecl^F$. Let $K$ be a countable subfield of $F$, and $C \subseteq F$ an $\ecl^F$-closed exponential subfield of $F$ containing $K$.

Then $C^K$ is a full, $\pcl^F_K$-closed $K$-powered subfield of $F^K$.
\end{prop}

\begin{proof}
Since $C$ is $\ecl^F$-closed, it inherits algebraic closedness and surjectivity of $\exp$ from $F$. Hence $C^K$ is a full $K$-powered field.

Now let $z \in F^n$. Then $\delta^K(z/C) \geqslant \delta^\exp(z/C) \geqslant 0$, with equality holding if and only if $z \in C^n$. Hence $C$ is $\pcl_K^F$-closed.
\end{proof}

\begin{coro}\label{pclccp}
Let $F$ be an exponential field with the countable closure property for $\ecl$, $K \subseteq F$ a countable subfield. Then $(F^K, \pcl^F_K)$ satisfies the countable closure property.

In particular, if $K$ is a countable subfield of $\mathbb{B}_\exp$ (resp$.$ $\mathbb{C}$) then $\mathbb{B}^K$ (resp$.$ $\mathbb{C}^K$) has the countable closure property.
\end{coro}

\begin{proof}
Let $A \subseteq F$ be  finite. Then $\pcl_K^F(A)\subseteq \ecl^F(A \cup K)$, which is countable.
\end{proof}

\section{Classification of Extensions}\label{classextsec}

In this section we introduce \textit{good bases} of finitely generated extensions, that is, bases which determine the extension up to isomorphism.

\begin{defn}
Let $D_1 \leq D_2$ be a finitely generated, kernel-preserving extension of partial $K$-powered fields. A \textit{good basis} for the extension, or a \textit{good basis of $D_2$ over $D_1$}, is a $\mathbb{Q}$-vector space basis $b$ of $D_2$ over $D_1$ such that if $b'$ is a tuple in some extension $D$ of $D_1$ satisfying $\Loc(b/D_1)=\Loc(b'/D_1)$, the partial $K$-powered field $\langle D_1, b' \rangle_\mathbb{Q}$ is isomorphic to $D_2$ over $D_1$.
\end{defn}

For any $k \in \mathbb{N}^+$, let $[k]:(F^\times)^n \rightarrow (F^\times)^n$ denote the componentwise multiplication-by-$k$ map.

\begin{defn}
Let $F$ be a field, $W$ an algebraic subvariety of $(F^\times)^n$ defined over $F$. We say $W$ is \textit{Kummer-generic} over $F$ if for all $k \in \mathbb{N}^+$, $[k]^{-1}(W)$ is irreducible over $F$.
\end{defn}

The notion of Kummer-generic is due to Hils \autocite[Definition 4.1]{Hil}, although he only defines Kummer-genericity over algebraically closed fields.

\begin{lem}\label{alggb}
    Let $D_1 \leq D_2$ be a finitely generated, kernel-preserving extension of partial $K$-powered fields, $b$ a basis for the extension. We have that $b$ is a good basis if and only if $\algloc \left( \exp\left(  b \right)/F_1 \right)$ is Kummer-generic over $F_1$.
\end{lem}

\begin{proof}
Let $b$ be a basis for $D_2$ over $D_1$, and $D$ an extension of $D_1$ with a basis $b'$ such that $\Loc(b/D_1)=\Loc(b'/D_1)$; we denote this common $K$-powers locus by $(L,W)$.

For any $k \in \mathbb{N}^+$, we have that $\algloc\left( \exp\left( \frac bk \right) /F_1 \right)$ is an irreducible component of $[k]^{-1}(W)$, and thus the $K$-powers loci of the points $\frac bk$ and $\frac {b'}k$ coincide for every $k$ if and only if $[k]^{-1}(W)$ has exactly one irreducible component over $F_1$ for each $k$, that is, if $W$ is Kummer-generic over $F_1$. Hence, we conclude by Lemma \ref{fingeniso}.
\end{proof}

Using Lemma \ref{alggb} we will prove that if $D$ is a $K$-powered field that is finitely generated or full then every finitely generated extension of $D$ has a good basis.

If $D$ is full, we get it as a consequence of the following result of Zilber.

\begin{thm}[{\autocite[Corollary 1.5]{Zil06}}, see also \cite{BGH}]\label{kummerthm}
Let $F$ be an algebraically closed field, $W \subseteq (F^\times)^n$. If $W$ is not contained in a coset of an algebraic subgroup of $(F^\times)^n$, then there is $m \in \mathbb{N}^+$ such that $[m]^{-1}(W)$ is Kummer-generic over $F$.
\end{thm}

If $D$ is finitely generated, then we follow the method of proof of \autocite[Section 3.3]{BK18}. Recall the notion of a \textit{division sequence} for a point in the multiplicative group of a field.

\begin{defn}
Let $F$ be a field, $w \in (F^\times)^n$. A \textit{division sequence} for $w$ is a sequence $(w_k)_{k \in \mathbb{N}^+}$ such that $w_1=w$ and for all $h,k \in \mathbb{N}^+$ we have $(w_{hk})^h=w_k$.

We denote by $\widehat{T}$ the group of division sequences of $1 \in (\mathbb{Q}^\textnormal{alg})^n$ with componentwise multiplication.
\end{defn}

\begin{prop}\label{otherkummer}
Let $D$ be a partial $K$-powered field with cyclic kernel that is either finitely generated or a finitely generated extension of a full $K$-powered field, let $D_1$ be a finitely generated, kernel-preserving extension of $D$, and let $b$ be a basis for the extension. 

Then there is $m \in \mathbb{N}^+$ such that $\algloc\left(\exp \left(\frac b{m}\right)/F\right)$ is Kummer-generic over $F$.  
\end{prop}

\begin{proof}
    Let $\xi_b:\textnormal{Gal}(F_1^\textnormal{alg}/F_1) \rightarrow \widehat{T}$ denote the Kummer map defined by $$\sigma \mapsto \left( \frac{\sigma \left( \exp \left(\frac{b}{k} \right) \right)}  {\exp \left(\frac bk \right)} \right)_{k \in \mathbb{N}^+}.$$ 

    The image of $\xi_b$ has finite index in $\widehat{T}$: this is proved in \autocite[Proposition A.9]{BHP} in the case of a finitely generated $K$-powered field and in \autocite[Section 3, Claim 2]{BGH} in the case of full $K$-powered field (see also \autocite[Proposition 3.24]{BK18}). Let $m$ be the exponent of the finite quotient $\widehat{T}/\im(\xi_b)$: then $\widehat{T}^m=\im(\xi_b)$.

    Let then $\xi_{\frac bm}$ be defined analogously to $\xi_b$ on $\textnormal{Gal}\left(F_1^\textnormal{alg}/F \left( \exp\left( \frac bm\right)\right) \right)$. We claim $\xi_{\frac bm}$ is surjective. In fact, if $t=(t_k)_{k \in \mathbb{N}^+}$ is a division sequence of 1, then $t^m \in \im(\xi_b)$ and thus there is $\sigma$ such that $\xi_b(\sigma)=t^m$. Hence we see that $\sigma\left( \exp \left(\frac bm\right)  \right)=\exp\left( \frac bm\right)$, that is, $\sigma \in \textnormal{Gal}\left(F_1^\textnormal{alg}/F\left(\exp\left(\frac bm\right)\right)\right)$, and that for every $k \in \mathbb{N}^+$ we have $$\sigma \left( \exp \left( \frac b{mk} \right) \right)=t_{mk}^m \exp \left( \frac b{mk} \right) = t_k \exp \left( \frac b {mk} \right)$$ where the first inequality holds by definition of $\xi_b$ and the second one by definition of division sequence. Hence, $\xi_{\frac bm}(\sigma)=t$. Since $t$ was arbitrary, this proves surjectivity.

    Then for all $k \in \mathbb{N}^+$, all the values of $[k]^{-1}\left( \exp \left( \frac bm \right) \right)$ are in the same orbit under $\textnormal{Gal}(F_1^\alg/F)$. Hence, $[k]^{-1}\left(\algloc \left( \exp \left( \frac bm \right) \right)\right)$ is irreducible over $F$.
\end{proof}

Putting these results together we obtain the existence of good bases.

\begin{prop}\label{goodbasesprop}
Let $D$ be a partial $K$-powered field with cyclic kernel that is either finitely generated or full. Every finitely generated, kernel-preserving extension of $D$ has a good basis.
\end{prop}

\begin{proof}
Let $b$ be a basis an extension $D \leq D_1$, and $(L,W) := \Loc(z/D_1)$. Since the extension is kernel-preserving, $\exp(b)$ does not satisfy any multiplicative relation, and as it is a generic point of $W$, this implies that $W$ is not contained in any coset of an algebraic subgroup.

If $D$ is full then $F$ is algebraically closed and by Theorem \ref{kummerthm}, there is some $m$ such that $[m]^{-1}(W)$ is Kummer-generic. If $D$ is finitely generated, then the same follows from Proposition \ref{otherkummer}. Either way, we conclude applying Lemma \ref{alggb}.
\end{proof}

\begin{coro}\label{coumanyext}
Let $D$ be a countable partial $K$-powered field with cyclic kernel that is either finitely generated or full. Up to isomorphism, there are only countably many finitely generated extensions of $D$.
\end{coro}

\begin{proof}
A finitely generated extension is determined up to isomorphism by the locus of a good basis. Since $D$ is countable, there are only countably many possible loci of points over it; therefore, there can only be countably many finitely generated extensions.
\end{proof}

We conclude by using good bases to prove two statements which will be needed later on.

\begin{lem}[Uniqueness of the generic type]\label{4.13}
Let $D_0 \treq D$ be an extension of $K$-powered fields, with $D_0$ partial and $D$ full. If $v \in D \setminus \pcl(D_0)$, then $\langle D_0, v \rangle_\mathbb{Q} \treq D$. Moreover, $v$ is a good basis for $\langle D_0, v \rangle_\mathbb{Q}$ over $D_0$, and hence $\langle D_0, v \rangle_{\mathbb{Q}} \cong \langle D_0, v' \rangle_{\mathbb{Q}}$ for all $v' \in D \setminus \pcl(D_0)$.
\end{lem}

\begin{proof}
Let $z \in D^n$. Then $\delta(z/D_0v)=\delta(vz/D_0)-\delta(v/D_0)=\delta(vz/D_0)-1$. By assumption $vz \notin (\pcl(D_0))^{n+1}$ and therefore $\delta(vz/D_0)>0$. Hence, $\delta(z/D_0v) \geqslant 0$, and $\langle D_0,v \rangle_\mathbb{Q} \treq D$.

To see that $v$ is a good basis we just notice that since it is not in $\pcl(D_0)$ its locus over $A$ must be the pair $(V,F^\times)$. $F^\times$ coincides with its preimages under multiplication-by-$k$ maps, and it is irreducible: hence $v$ is a good basis.
\end{proof}

\begin{prop}[Existence and uniqueness of full closures]\label{fullclosures}
Let $D$ be a partial $K$-powered field with cyclic kernel. Then there exists a full $K$-powered field $D^f$ that strongly extends $D$, such that there are no intermediate full $K$-powered fields $D_1$ with $D \subseteq D_1 \subseteq D^f$.

Moreover, if $D$ is finitely generated, or a finitely generated extension of a countable full $K$-powered field, then $D^f$ is unique up to isomorphism.
\end{prop}

\begin{proof}
We only sketch the proof, as it is essentially the same as \autocite[Theorem 2.18]{Kir13} (see also \autocite[Theorem 4.17]{BK18}.)

For existence, we embed $V$ into a large $K$-vector space $\mathcal{V}$ and $F$ in a large algebraically closed field $\mathcal{F}$. We extend the exponential with an iterated procedure: we coherently map the $\mathbb Q$-span of any element of $V$ that is not in $D$ to rational powers of an element of $\mathcal{F}$ that is transcendental over $F$. Moreover, for each element $w \in \mathcal{F}^\times$ that is algebraic over the image of $\exp$ at any given step, we take an element in $\mathcal{V}$ that is not in $V$ and map its $\mathbb Q$-span (again, coherently) to rational powers of $w$. This procedure yields a strong extension, and iterating we produce a full $K$-powered field. It is clear from the construction that no intermediate $D_1$ exists.

If $D$ is finitely generated or a finitely generated extension of a countable full $K$-powered field, we assume that there are two full $K$-powered field containing $D$ as in the statement. We then break up one of them into a countable union of finitely generated extensions of $D$, and embed each of these extensions into the other full $K$-powered field using the good bases given by Proposition \ref{goodbasesprop}. These embeddings will produce an isomorphism.
\end{proof}

Given a partial $K$-powered field $D$, any full $K$-powered field obtained as in the proof of Proposition \ref{fullclosures} will be called a \textit{full closure} of $D$.

\begin{lem}\label{stronginfull}
    Let $D_0$ be a full $K$-powered field, $D$ a purely powers-transcendental extension of $D_0$.

    If $D^f$ is a full closure of $D$, then $D_0 \treqcl D^f$.
\end{lem}

\begin{proof}
    $D^f$ is obtained from $D$ by an iterated construction, so it it sufficient to check that at every step we have a purely powers-transcendental extension of $D_0$. If $z \in V \setminus D$ then $\delta(z/D_0)=1$. Similarly, if $w \in F^\times \setminus \exp(D)$ and $z' \in D^f$ is a point with $\exp(z')=w$, then $\delta(z'/D_0)=1$. Hence the statement holds.
\end{proof}

\section{Amalgamation and Excellence}\label{a and e section}

Fra\"iss\'e's original amalgamation theorem \autocite[Chapter 7]{Hod} gives the existence and uniqueness of a countable $\aleph_0$-saturated model $\mathbb{F}$ in a class $\mathcal{C}$ of relational structures which is specified by saying that all finite substructures are from a certain subclass. The uniqueness is by the older back-and-forth technique, and the existence is from an amalgamation construction.

The first aim of this section is to use a variant of Fra\"iss\'e's Theorem to show that two categories of $K$-powered fields, denoted $\mathcal{C}(D_0)$ and $\mathcal{C}^\tr(D_0)$, have Fra\"iss\'e limits $\mathbb{F}^K(D_0)$ and $\mathbb{F}^{K,\tr}(D_0)$.

The second aim is to construct quasiminimal $K$-powered fields $\mathbb{E}^K(D_0)$ and $\mathbb{E}^{K,\tr}(D_0)$ of size continuum, and the last aim is to characterize these models semantically in the categories $\mathcal{C}(D_0)$ and $\mathcal{C}^\tr(D_0)$. We achieve these using the theory of \textit{quasiminimal excellent classes} developed in \cite{Zil05}, \cite{Kir10b}, \cite{BHHKK}, showing that the Fra\"iss\'e limits $\mathbb{F}^K(D_0)$ and $\mathbb{F}^{K,\tr}(D_0)$ generate such classes. 

\subsection{Amalgamation to a Fra\"iss\'e limit}

There are many variants and generalisations of Fra\"iss\'e's theorem, all with essentially the same proof. We use a version from \autocite[Theorem 2.18]{Kir09}, extending \cite{DG92}, which allows us to replace the class $\mathcal{C}$ by a category of countable structures where we can specify which embeddings we allow, and finite becomes finitely generated in a suitable sense.

\begin{defn}\label{category}
Let $K$ be a countable field, $D_0$ be a partial $K$-powered field with cyclic kernel.

\begin{itemize}
    \item[1.] $\mathcal{C}(D_0)$ is the category of partial $K$-powered fields strongly extending $D_0$, with strong embeddings of $D_0$ as the morphisms;
    \item[2.] $\mathcal{C}^\tr(D_0)$ is the full subcategory of $\mathcal{C}(D_0)$ whose objects are purely powers-transcendental extensions of $D_0$.
\end{itemize}
\end{defn}

\begin{defn}
Let $D_0$ be a partial $K$-powered field, and let $\mathcal{C}$ be one of the categories $\mathcal{C}(D_0)$ and $\mathcal{C}^\tr(D_0)$. An object $A$ of $\mathcal{C}$ is \textit{$\aleph_0$-saturated in $\mathcal{C}$} if for every finitely generated subobject $A_1$ and every finitely generated extension $A_2$ of $A_1$ in $\mathcal{C}$, $A_2$ embeds into $A$ over $A_1$ in $\mathcal{C}$.

A \textit{Fra\"iss\'e limit} for $\mathcal{C}$ is a countable, $\aleph_0$-saturated object of $\mathcal{C}$.
\end{defn}

\begin{rem}
 The usual back-and-forth argument shows that if an object $D$ is $\aleph_0$-saturated in $\mathcal{C}$ then it is also $\aleph_1$-universal in $\mathcal{C}$ (every object generated by a set of cardinality at most $\aleph_0$ embeds in $M$) and $\aleph_0$-homogeneous in $\mathcal{C}$ (every isomorphism between finitely generated substructures extends to an automorphism).   
\end{rem}

\begin{thm}
Let $\mathcal{C}$ be $\mathcal{C}(D_0)$ for a finitely generated or countable full $K$-powered field $D_0$, or $\mathcal{C}^\tr(D_0)$ for a countable full $K$-powered field $D_0$. 

Then $\mathcal{C}$ has a Fra\"iss\'e limit.
\end{thm}

\begin{nota}
We denote by $\mathbb{F}(D_0)$ and $\mathbb{F}^\tr(D_0)$ the Fra\"iss\'e limits of $\mathcal{C}(D_0)$ and $\mathcal{C}^\tr(D_0)$ respectively. We write just $\mathbb{F}^K$ for $\mathbb{F}(SB^K)$.
\end{nota}

\begin{proof}
By \autocite[Theorem 2.18]{Kir09} it suffices to prove that $\mathcal{C}$ is an \textit{amalgamation category} in the sense of \autocite[Definition 2.17]{Kir09}. We use the numbering of the axioms of an amalgamation category from \autocite[Definition 5.3]{BK18}.

(AC1) translates here as all morphisms in the category being embeddings and (AC2) to the category being closed under unions of $\omega$-chains; these are verified immediately. 

(AC4) states in this context that every finitely generated object of $\mathcal{C}$ has countably many finitely generated extensions up to isomorphism, and it holds by Corollary \ref{coumanyext}. (AC3) states that there are only countably many finitely generated objects in $\mathcal{C}$ up to isomorphism, and it follows from (AC4) as $D_0$ is an object of $\mathcal{C}$.

(AC5) is the amalgamation property for finitely generated objects, and it implies (AC6), the joint embedding property for finitely generated objects. Hence it remains to prove (AC5).

Let $D_1,D_2,D_3$ be finitely generated objects of $\mathcal{C}$, with strong embeddings of $D_1$ into $D_2$ and $D_3$. By Proposition \ref{fullclosures}, each of the objects has a full closure that is unique up to isomorphism, which we denote by $D_1^f,D_2^f,D_3^f$.

By Proposition \ref{fullclosures}, the full closures of $D_1$ inside $D_2^f$ and $D_3^f$ are both isomorphic to $D_1^f$, so we can choose strong embeddings as in the dashed arrows in the diagram. Thus we let $D_4$ be the free amalgam of $D_2^f$ and $D_3^f$ over $D_1^f$.

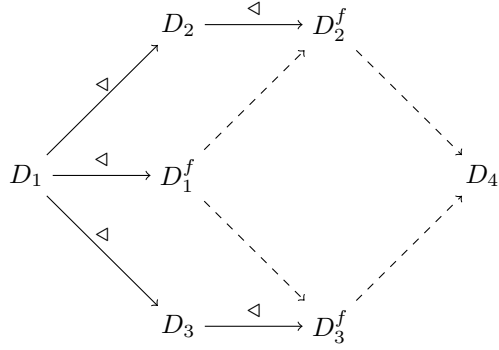
$$\begin{tikzpicture} [node distance=2cm]
\node (a) at (-2,0) {$D_1$};
\node (b) at (0,0) {$D_1^f$};
\node (c) at (4,0) {$D_4$};
\node (d) at (0,2) {$D_2$};
\node (e) at (0,-2) {$D_3$};
\node (f) at (2,2) {$D_2^f$};
\node (g) at (2,-2) {$D_3^f$};

\draw [->] (a) -- (b) node[midway, above]{$\treq$};
\draw [->] (a) -- (d) node[midway, above]{$\treq$};
\draw [->] (a) -- (e) node[midway, above]{$\treq$};
\draw [dashed, ->] (b) -- (f);
\draw [dashed, ->] (b) -- (g);
\draw [dashed, ->] (f) -- (c);
\draw [dashed, ->] (g) -- (c);
\draw [->] (e) -- (g) node[midway, above]{$\treq$};
\draw [->] (d) -- (f) node[midway, above]{$\treq$};
\end{tikzpicture}$$

By Proposition \ref{fullamalgam}, $D_4$ is a strong extension of $D_2^f$ and $D_3^f$, and thus it is a strong extension of $D_2$ and $D_3$. We take $D=\hull{D_2 \cup D_3}^{D_4}_K$, which is a finitely generated strong extension of $D_2$ and $D_3$ over $D_0$, and is the amalgam we need. 

Now assume $D_1, D_2$ and $D_3$ are purely powers-transcendental extensions of the full $K$-powered field $D_0$. By Lemma \ref{stronginfull} $D_1^f, D_2^f$ and $D_3^f$ are also purely powers-transcendental over $D_0$. Then by Proposition \ref{pptamal}, $D_4$ is a purely powers-transcendental extension of $D_0$. Hence $D$ is a finitely generated, purely powers-transcendental extension of $D_0$, with embeddings as in the definition of the amalgamation property.
\end{proof}

\subsection{The models $\mathbb{E}^K(D_0)$ and $\mathbb{E}^{K,\tr}(D_0)$}

\begin{defn}
Let $M$ be an $L$-structure for a countable language $L$, equipped with a pregeometry $\cl$. We say that $M$ is a \textit{quasiminimal pregeometry structure} if it satisfies the following axioms, in which $\qftp$ denotes the quantifier-free type:

\begin{itemize}
\item[(QM1)] (The pregeometry is determined by the language) If $a$ and $b$ are finite tuples with $\qftp(a)=\qftp(b)$, then $a$ and $b$ have the same $\cl$-dimension.
\item[(QM2)] (Infinite dimensionality) $M$ is infinite-dimensional with respect to $\cl$.
\item[(QM3)] (Countable closure property) If $A \subseteq M$ is finite, then $\cl(A)$ is countable.
\item[(QM4)] (Uniqueness of the generic type) If $C,C'$ are countable closed substructures of $M$, they are enumerated so that $\qftp(C)=\qftp(C')$, $a \in M \setminus C$ and $a' \in M' \setminus C'$, then $\qftp(C,a)=\qftp(C',a')$.
\item[(QM5a)] ($\aleph_0$-homogeneity over the empty set) Let $a,b$ be finite tuples in $M$ such that $\qftp(a)=\qftp(b)$, and let $a' \in \cl(a)$. Then there is $b' \in \cl(b)$ such that $\qftp(a,a')=\qftp(b',b')$.
\item[(QM5b)] (Non-splitting over a finite set) Let $C$ be a closed subset, and let $b$ be a finite tuple. Then there is a finite tuple $c$ in $C$ such that for all finite tuples $a,a'$ in $C$, if $\qftp(a/c)=\qftp(a'/c)$ then $\qftp(a/cb)=\qftp(a'/cb)$.
\end{itemize}

A \textit{weakly quasiminimal pregeometry structure} is an $L$-structure which satisfies all the axioms except possibly (QM2) (it is not necessarily infinite dimensional).
\end{defn}

\begin{defn}
Given weakly quasiminimal pregeometry structures $(M_1,\cl^{M_1})$ and $(M_2,\cl^{M_2})$ (in the same language), an embedding $\theta:M_1\hookrightarrow M_2$ is a \textit{closed embedding} if for every $X \subseteq M_1$, $\cl^{M_2}(\theta(X))=\theta(\cl^{M_1}(X))$.

Given a quasiminimal pregeometry structure $M$ in a language $L$, we denote by $\mathcal{K}(M)$ the smallest class of $L$-structures which contains $M$ and all its closed substructures and is closed under isomorphisms and under taking unions of directed systems of closed embeddings. Such a class is called the \textit{quasiminimal class} attached to $M$.
\end{defn}

\begin{fact}[{\autocite[Theorem 2.3]{BHHKK}}]\label{categoricity}
If $\mathcal{K}$ is a quasiminimal class, then every structure $A \in \mathcal{K}$ is a weakly quasiminimal pregeometry structure. For every cardinal $\kappa$, there is exactly one model of dimension $\kappa$ in $\mathcal{K}$ (up to isomorphism); in particular, $\mathcal{K}$ is uncountably categorical.
\end{fact}

\begin{defn}
Let $K$ be a countable field, and let $\mathcal{L}^K$ be the language of $K$-powered fields from Definition \ref{language}.

Let $D_0$ be a countable $K$-powered field that is either finitely generated or full.

We extend this to a language $\mathcal{L}^{K,QE}(D_0)$ by adding constant symbols for all elements of $D_0$ and, for every pair $(L,W)$ where $L$ is a $K$-linear subspace of $V_0^{n+k}$ and $W$ is an algebraic subvariety of $F_0^{n+k}$ defined over $\mathbb{Q}$, a $k$-ary predicate $\phi_{L,W,n}(x)$. We interpret $\phi_{L,W,n}(a)$ as $$\exists b \left[ (b,a,\exp(b,a)) \in L \times W \wedge \ldim_\mathbb{Q}(b/D_0a)=n\right].$$
\end{defn}

\begin{lem}\label{hullqftp}
    Let $D$ be an object of $\mathcal{C}(D_0)$ or $\mathcal{C}^\tr(D_0)$, and let $a,b \in D^k$. Then $\qftp(a)=\qftp(b)$ if and only if there is an isomorphism $\theta:\hull{D_0,a}^D_K \cong \hull{D_0,b}^D_K$ fixing $D_0$ pointwise and sending $a$ to $b$.
\end{lem}

\begin{proof}
$(\Leftarrow)$ is obvious.

For $(\Rightarrow)$, let $a,b$ be $k$-tuples in $\mathbb{F}$ with the same $\mathcal{L}^{K,QE}(D_0)$-quantifier free type. Let $D_a=\hull{D_0,a}^\mathbb{F}_K$, and  let $a' \in D_a^n$ be a good basis for $D_a$ over $D_0$. Then let $(L,W):=\Loc(a',a/D_0)$. By definition of $\mathcal{L}^{K,QE}(D_0)$, $\mathbb{F} \vDash \phi_{L,W,n}(a)$ and therefore since $a$ and $b$ have the same quantifier-free type, $\mathbb{F} \vDash \phi_{L,W,n}(b)$, that is, there is $b' \in \mathbb{F}^n$ such that $(b',b,\exp(b'),\exp(b)) \in L \times W$ and $\ldim_{\mathbb{Q}}(b'/D_0b)=n$. We want to show that $(L,W)=\Loc(b',b/D_0)$, so let $(L', W')$ be the $K$-powers locus of $(b',b)$ over $D_0$. Again, $a$ and $b$ have the same quantifier-free type so there is $a'' \in \mathbb{F}^n$ such that $(a'',a, \exp(a''), \exp(a)) \in L' \times W'$ and $\ldim_{\mathbb{Q}}(a''/D_0a)=n$. Then
\begin{align*}
	\dim L' + \dim W' & = \delta(a'',a/D_0) + \ldim_\mathbb{Q}(a'',a/D_0) \\
	                  & = \delta(a''/D_0) + n + \ldim_\mathbb{Q}(a/D_0) \\
	                  & \geqslant \delta (a'/D_0) + \ldim_{\mathbb{Q}}(a',a/D_0) \\
	                  & =\delta(a',a/D_0) + \ldim_{\mathbb{Q}}(a',a/D_0) \\
	                  & = \dim L + \dim W.
\end{align*} 

Therefore, we have that $(L', W')=(L,W)=\Loc(b',b/D_0)$, and from this we deduce that $\Loc(a'/D_0)=\Loc(b'/D_0)$. Since $a'$ is a good basis for $D_a$, there is an isomorphism between $D_a$ and the $K$-powered field $\langle D_0, b' \rangle_{\mathbb Q}$ which fixes $D_0$ and maps $a'$ to $b'$ and hence $a$ to $b$. Note that by the above $\delta(b',b/D_0)=\delta(a',a/D_0)$, and by the same proof for every tuple $b''$ in $D$ we must have $\delta(b'',b/D_0) \geqslant \delta(b',b/D_0)$. Hence $\langle D_0, b \rangle_\mathbb{Q} =\hull{D_0, b}_K^D$.
\end{proof}

\begin{thm}\label{qps}
Let $K$ be a countable field; let $\mathbb{F}$ be one of the Fra\"iss\'e limits $\mathbb{F}^K(D_0)$ and $\mathbb{F}^{K,\tr}(D_0)$ (the latter only if $D_0$ is a countable full $K$-powered field).

Then $\mathbb{F}$ in the language $\mathcal{L}^{K,QE}(D_0)$ and with the pregeometry $\pcl_K^{\mathbb{F}}$ is a quasiminimal pregeometry structure.
\end{thm}

\begin{proof}
We denote by $\mathcal{C}$ the category of which $\mathbb{F}$ is the Fra\"iss\'e limit (so $\mathcal{C}(D_0)$ or $\mathcal{C}^\tr(D_0)$).

The arguments are then similar to the proof of \autocite[Theorem 6.9]{BK18}. 

(QM3) holds because $\mathbb{F}$ is countable.

For (QM1) and (QM5a) we use Lemma \ref{hullqftp}: let $a$ and $b$ be tuples in $\mathbb{F}$ with $\qftp(a)=\qftp(b)$, $D_a=\hull{D_0,a}^{\mathbb{F}}_K$ and $D_b=\hull{D_0,a}^{\mathbb{F}}_K$. Then there is an isomorphism $\theta:D_a \rightarrow D_b$ which fixes $D_0$ pointwise and maps $a$ to $b$, hence $\pdim^{D_a}(a/D_0)=\pdim^{D_b}(b/D_0)$. By Lemma \ref{strong ext preserve ktd} then $\pdim^{\mathbb{F}}(a/D_0)=\pdim^{\mathbb{F}}(b/D_0)$, so (QM1) holds. For (QM5a) we have that, since $\mathbb{F}$ is a Fra\"iss\'e limit and hence is $\aleph_0$-homogeneous, $\theta$ extends to an automorphism $\theta'$ of $\mathbb{F}$, so for every $a' \in \pcl_K^{\mathbb{F}}(a)$ there is $b' \in \pcl_K^{\mathbb{F}}(b)$ such that $\qftp(a,a')=\qftp(b,b')$.

Now we prove (QM2). Let $n \in \mathbb{N}$. Then there is a strong extension $D_n$ of $D_0$ which is generated by an $n$-tuple $a$ that is generic over $D_0$, so $\pdim^{D_n}(a/D_0)=\delta^K(a/D_0)=n$; $D_n$ is an object of $\mathcal{C}$. Since $\mathbb{F}$ is $\aleph_0$-universal in $\mathcal{C}$, $D_n$ embeds into $\mathbb{F}$, which thus contains a tuple $a$ with $\pdim^\mathbb{F}(a/D_0)=n$. Hence $\mathbb{F}$ is infinite dimensional.

For (QM4), suppose $D_1$ and $D_2$ are countable $\pcl$-closed $K$-powered subfields of $\mathbb{F}$, enumerated so that $\qftp(D_1)=\qftp(D_2)$, and let $a \in \mathbb{F} \setminus D_1$ and $b \in \mathbb{F} \setminus D_2$. By Lemma \ref{4.13}, the extensions $D_1 \treq \langle D_1, a \rangle_\mathbb Q$ and $D_2 \treq \langle D_2, b \rangle_\mathbb Q$ are isomorphic, hence $\qftp(D_1, a)=\qftp(D_2,b)$.

Finally we prove (QM5b). Suppose $D$ is a $\pcl$-closed $K$-powered subfield of $\mathbb{F}$, and let $b$ be a finite tuple in $\mathbb{F}$; without loss of generality we assume that $b$ is a good basis for the extension $D \treq \hull{D,b}^\mathbb{F}_K =:D_b$. Let $c$ be a finite tuple in $D$ such that $\Loc(b/D)$ is defined over $\langle D_0, c \rangle_\mathbb Q$, and let $a,a'$ be finite tuples in $D$ such that $\qftp(a/c)=\qftp(a'/c)$. By Lemma \ref{hullqftp}, there is an isomorphism $\theta_0:D_{a,c} \cong D_{a',c}$, where $D_{a,c}=\hull{D_0,a,c}^\mathbb{F}_K$ and $D_{a',c}=\hull{D_0,a',c}^\mathbb{F}_K$. Replacing $b$ by one of its rational multiples if necessary, so that it is a good basis of $D_{a,c} \leq \langle D_{a,c}, b \rangle_{\mathbb{Q}}$, we have that $\theta_0$ extends to an isomorphism $\theta_1: \langle D_{a,c}, b \rangle_\mathbb Q \cong \langle D_{a',c}, b \rangle_\mathbb Q$.

Now we claim that $\langle D_{a,c}, b \rangle_\mathbb Q \treq D_b$. To see this, let $z$ be a finite tuple in $D_b$: as we assumed that $b$ is a basis for the extension $D \treq D_b$, $z$ is the sum of a tuple $z_D$ in $D$ and a tuple $z_b$ in $\langle b \rangle_\mathbb Q$, and thus $\delta(z/D_{a,c}, b)=\delta(z_D/D_{a',c},b)$. However, $z_D$ and $D_{a,c}$ are contained in the $\pcl$-closed $K$-powered field $D$, while $b$ is $\mathbb{Q}$-linearly independent over $D$, so $\delta(z_D/D_{a,c}, b)=\delta(z_D/D_{a,c})$. Since $D_{a,c} \treq D$ this is non-negative. Hence $\langle D_{a,c}, b \rangle_\mathbb Q \treq \langle D, b \rangle_\mathbb{Q}$ as we wanted. By the same argument, $\langle D_{a',c},b \rangle_\mathbb Q$ is also strong in $D_b$, and since the latter is strong in $\mathbb{F}$ we have that  $\langle D_{a,c}, b \rangle_\mathbb Q$  and $\langle D_{a',c}, b \rangle_\mathbb Q$ are finitely generated objects in $\mathcal{C}$, strong in $\mathbb{F}$, with the isomorphism $\theta_1$ between them. Hence there is an isomorphism $\hull{D_0,a,b,c}_K^\mathbb{F} \cong \hull{D_0,a',b,c}_K^\mathbb{F}$, which implies $\qftp(a'/bc)=\qftp(a/bc)$, as we wanted. 
\end{proof}

By Fact \ref{categoricity} and Theorem \ref{qps}, each of the quasiminimal classes generated by $\mathbb{F}^K(D_0)$ and $\mathbb{F}^{K,\tr}(D_0)$ contains a unique model of cardinality continuum (up to isomorphism). We denote these by $\mathbb{E}^K(D_0)$ and $\mathbb{E}^{K,\tr}(D_0)$. As with the Fra\"iss\'e limit, we write $\mathbb{E}^K$ for $\mathbb{E}^K(SB^K)$.

\subsection{Algebraic saturation}

In this subsection we introduce a notion of \textit{algebraic saturation} and use it to characterize the models $\mathbb{E}^K(D_0)$ and $\mathbb{E}^{K,\tr}(D_0)$.

\begin{defn}
    Let $\mathcal{C}$ denote $\mathcal{C}(D_0)$ or $\mathcal{C}^\tr(D_0)$. An object $D$ in $\mathcal{C}$ is \textit{algebraically saturated in $\mathcal{C}$} if for all finitely generated $D_1 \treq D$, and all finitely generated, powers-algebraic extensions $D_1 \treq D_2$ in $\mathcal{C}$, $D_2$ embeds (strongly) in $D$ over $D_1$.
\end{defn}

\begin{thm}\label{characterize}
    Let $K$ be a countable field. Let $D_0$ be a finitely generated or full, countable $K$-powered field (resp$.$ a full, countable $K$-powered field). 

    The $K$-powered field $\mathbb{E}^K(D_0)$ (resp$.$ $\mathbb{E}^{K,\tr}(D_0)$) is up to isomorphism the unique full $K$-powered field of cardinality continuum which strongly extends $D_0$ (resp$.$ is a purely powers-transcendental extension of $D_0$), is algebraically saturated in $\mathcal{C}(D_0)$ (resp. $\mathcal{C}^\tr(D_0)$) and has the countable closure property.
\end{thm}

\begin{proof}
Let $\mathbb{E}$ denote $\mathbb{E}^K(D_0)$ or $\mathbb{E}^{K,\tr}(D_0)$, and let $\mathcal{C}$ and $\mathbb{F}$ accordingly denote $\mathcal{C}(D_0)$ or $\mathcal{C}^\tr(D_0)$ and $\mathbb{F}^K(D_0)$ or $\mathbb{F}^{K,\tr}(D_0)$.

$(\Rightarrow)$ As $\mathbb{F}$ is $\aleph_0$-saturated, it is also algebraically saturated, and thus any of its $\pcl$-closed substructures is also algebraically saturated. Since $\mathbb{E}$ is the union of a directed system of closed embeddings of closed substructures of $\mathbb{F}$ the statement follows.

$(\Leftarrow)$ Let $D$ be a full object of $\mathcal{C}$ of size continuum that is algebraically saturated in $\mathcal{C}$ and has the countable closure property. Any finite $\pcl$-dimensional, $\pcl$-closed substructure $C$ of $D$ is countable, and hence it embeds into $\mathbb{F}$ as a $\pcl$-closed substructure by $\aleph_1$-universality of the Fra\"iss\'e limit. Hence $C$ lies in the quasiminimal class $\mathcal{K}(\mathbb{F})$ As $D$ is the union of its finite $\pcl$-dimensional, $\pcl$-closed substructures, $D$ is also in $\mathcal{K}(\mathbb{F})$ and is therefore isomorphic to $\mathbb{E}$.
\end{proof}

\section{$K$-powers-closed fields}\label{k-pow-closed fields sections}

In this section we introduce a notion of \textit{strong $K$-powers closedness}, prove that it is equivalent to algebraic saturation in the category $\mathcal{C}(D_0)$, and that when $D_0$ is finitely generated, so we have a Schanuel property in $\mathcal{C}(D_0)$, it can be reduced to the simpler notion of \textit{$K$-powers-closedness}. 

\subsection{Classification of strong extensions}

\begin{defn}
Let $D$ be a partial $K$-powered field, $L \leq D^n$ a $K$-linear subspace and $W \subseteq (F^\times)^n$ an algebraic subvariety. 

The pair $(L,W)$ is \textit{free} if $L$ is not contained in any $\mathbb{Q}$-affine subspace of $D^n$ and $W$ is not contained in a coset of a proper algebraic subgroup of $(F^\times)^n$. 

The pair $(L,W)$ is \textit{rotund} if for every $\mathbb{Q}$-linear subspace $Q$ of $D^n$ with projections $\pi_Q:D^n \twoheadrightarrow D^n/Q$ and $\pi_{\exp(Q)}:(F^\times)^n \twoheadrightarrow (F^\times)^n/(\exp(Q))$ we have $$\dim \pi_Q(L) + \dim \pi_{\exp(Q)} (W) \geqslant n - \dim Q$$ and it is \textit{strongly rotund} if for every proper $Q$ we have $$\dim \pi_Q(L) + \dim \pi_{\exp(Q)} (W) > n - \dim Q.$$
\end{defn}

\begin{rem}
By taking $Q=\langle 0 \rangle_\mathbb{Q}$ we have that in particular a rotund pair always satisfies $\dim L + \dim W \geqslant n$ and a strongly rotund pair always satisfies $\dim L + \dim W > n$.
\end{rem}

\begin{prop}\label{classify}
Let $D$ be a partial $K$-powered field, $L \leq V^n$ a $K$-linear subspace, $W \subseteq (F^\times)^n$ an algebraic subvariety.

If the pair $(L,W)$ is free, there is a finitely generated, kernel-preserving extension $D_1$ of $D$, with a basis $b$ of $D_1$ over $D$ such that $\Loc(b/V)=(L,W)$. Moreover:

\begin{itemize}
    \item[$(1)$] The extension is strong if and only if $(L,W)$ is rotund;
    \item[$(2)$] If the extension is strong, it is powers-algebraic if and only if $\dim L + \dim W =n$;
    \item[$(3)$] The extension is purely powers-transcendental if and only if $(L,W)$ is strongly rotund.
\end{itemize}
\end{prop}

\begin{proof}
The first part of the statement is by the same argument as \autocite[Lemma 3.4]{Zil03}.

For the ``moreover'' part, (1) and (3) are obtained by the same argument as \autocite[Proposition 7.3]{BK18}; (2) follows easily.
\end{proof}

\subsection{$K$-powers-closedness}

\begin{defn}
Let $D_0 \treq D$ be an extension of partial $K$-powered fields.

$D$ is \textit{$K$-powers-closed} if for every free, rotund, $n$-dimensional pair $(L,W)$ in $D^n$, we have that $\exp(L) \cap W$ is Zariski-dense in $W$.

$D$ is \textit{strongly $K$-powers-closed over $D_0$} if for every free, rotund, $n$-dimensional pair $(L,W)$ in $D^n$ and every finite tuple $a$ in $D$ there is $b \in D^n$ such that $(b,\exp(b)) \in L \times W$ and $\ldim_\mathbb{Q}(b/D_0a)=n$.
\end{defn}

The main result obtained by the first author in \cite{Gal22} is that $\mathbb{C}^\mathbb{C}$ is $\mathbb{C}$-powers-closed.

\begin{prop}\label{saturation and pow closedness}
Let $D_0$ be either a finitely generated partial $K$-powered field or a full countable $K$-powered field. Recall that $\mathcal{C}(D_0)$ denotes the category of partial $K$-powered fields strongly extending $D_0$.

An object in $\mathcal{C}(D_0)$ is algebraically saturated in $\mathcal{C}(D_0)$ if and only if it is full and strongly $K$-powers-closed over $D_0$.
\end{prop}

\begin{proof}
    $(\Rightarrow)$ Let $D$ be an object of $\mathcal{C}(D_0)$ that is algebraically saturated in $\mathcal{C}(D_0)$. It is straightforward to show that $D$ is full, so we show that it is strongly $K$-powers-closed.
    
    Let $(L,W)$ be a free, rotund, $n$-dimensional pair in $D^n$, and $a \in D^k$ a tuple. Replacing $a$ by a basis for $\hull{D_0, a}_K^D$ if necessary, we may assume that $D_1:=\langle D_0, a \rangle_{\mathbb{Q}} \treq D$. Add to $D_1$ a point $b$ such that $(b,\exp(b))$ is generic in $L \times W$ over $D_1$; by Proposition \ref{classify}, this generates a strong, powers-algebraic extension $D_2$ of $D_1$. By algebraic saturation of $D$, $D_2$ embeds into $D$ over $D_1$. Hence there is a point $b'$, the image of $b$ under the embedding, such that $(b',\exp(b')) \in L \times W$, and (by genericity) $\ldim_\mathbb{Q}(b'/D_1)=n$. Hence $D$ is strongly $K$-powers-closed.

    $(\Leftarrow)$ Let $D$ be an object of $\mathcal{C}(D_0)$ that is full and strongly $K$-powers-closed over $D_0$. Let $D_1=\langle D_0, a\rangle_\mathbb{Q}$ be a finitely generated, strong extension of $D_0$ that is strong in $D$, and let $D_2$ be a finitely generated, strong, powers-algebraic extension of $D_1$. By Proposition \ref{goodbasesprop}, we may find a good basis $b$ for $D_2$ over $D_1$, so $D_2=\langle D_1, b \rangle_\mathbb{Q}$.

    We prove that $D_2$ strongly embeds into $D$ over $D_1$ by induction on $n:=\dim_{\mathbb{Q}}(D_2/D_1)$. Let $(L,W)=\Loc(b/D_1)$. 
    
    If $(L,W)$ is free, since the extension $D_1 \treq D_2$ is strong and powers-algebraic, by Proposition \ref{classify} we have that $(L,W)$ is rotund and $n$-dimensional. By strong $K$-powers-closedness of $D$ there is a point $b' \in D^n$ such that $(b', \exp(b')) \in L \times W$ and $\ldim_{\mathbb{Q}}(b'/D_1)=n$. Since $D$ strongly extends $D_1$, it must be the case that $\ldim_K(b'/D_1)+\trdeg(\exp(b')/F_1) \geqslant n$. As $\dim L + \dim W =n$, this implies that $(L,W)=\Loc(b'/D_1)$. Since $b$ is a good basis, the extension $\langle D_1, b' \rangle_\mathbb{Q}$ of $D_1$ is then isomorphic to $D_2$ over $D_1$, and it is contained in $D$, so $D_2$ embeds in $D$ over $D_1$.

    If $(L,W)$ is not free, then without loss of generality we may assume that either $L$ is contained in a $\mathbb{Q}$-affine subspace over $D_1$ of the form $z_i=c$ for some $c \in V_1$, or $W$ is contained in an algebraic subvariety of the form $w_i=d$ for some $d \in F_1^\alg$. In the first case we replace $D_1$ by $D_1':=\hull{D_1,c}_K^D$ and consider $D_2$ as a powers-algebraic strong extension of $D_1'$ of lower dimension (hence we may apply induction); the second case is treated analogously.
\end{proof}

It is clear that strong $K$-powers-closedness implies $K$-powers-closedness. The rest of this subsection will be dedicated to proving the converse under a transcendence assumption.

\begin{prop}\label{pow closed implies spow closed}
Let $D$ be a full $K$-powered field. If $D$ has a strong finitely generated partial $K$-powered subfield $D_0$ and it is $K$-powers-closed, then it is strongly $K$-powers-closed over $D_0$, and hence also algebraically saturated in $\mathcal{C}(D_0)$.
\end{prop}

We recall some basics from the theory of \textit{atypical intersections}.

\begin{defn}
    Let $W \subseteq (F^\times)^n$ be an algebraic subvariety, $J$ a coset of an algebraic subgroup of $(F^\times)^n$. An irreducible component $X$ of the intersection $W \cap J$ is \textit{typical} if $\dim X = \dim W + \dim J -n$, and it is \textit{atypical} otherwise.
\end{defn}

 We recall the following result, usually referred to as the \textit{weak CIT} (weak Conjecture on Intersections with Tori) or \textit{weak multiplicative Zilber-Pink}.

\begin{thm}[{\autocite[Corollary 3]{Zil02}}, see also {\autocite[Theorem 4.6]{Kir09}}]\label{weak cit}
    Let $W \subseteq (F^\times)^n$ be an algebraic subvariety. There is a finite set $\mathcal{T}_W=\{J_1,\dots,J_l\}$ of algebraic subgroups of $(F^\times)^n$ such that for any coset $c \cdot J$ of an algebraic subgroup $J$ of $(F^\times)^n$ and any atypical component $X$ of an intersection $W \cap c \cdot J$, there is $J_i \in \mathcal{T}_W$ such that $X$ is contained in a coset $w \cdot J_i$.

    Moreover, $X$ is typical with respect to $c \cdot J$, meaning that $$\dim X= \dim (W \cap  c \cdot J) + \dim ((w \cdot J_i) \cap (c \cdot J)) - \dim (c \cdot J).$$
\end{thm}

\begin{rem}[See {\autocite[Lemma 4.1]{Zil11}}]\label{w is algebraic}
  With notation as in the statement of Theorem \ref{weak cit}, if $W$ is defined over a field $F_0$ then we may assume that $w$ is algebraic over $F_0(c)$. This is because the irreducible component $X$ of $W \cap c \cdot J$ is an algebraic variety defined over $F_0(c)$, and therefore it contains points that are algebraic over this field.   
\end{rem}

\begin{lem}\label{first typical fibre lemma}
Let $F$ be a field, $W \subseteq (F^\times)^n$ an algebraic subvariety, $J \leq (F^\times)^n$ an algebraic subgroup. There is a Zariski-closed proper subset $W'$ of $W$ such that for each $w \in W \setminus W'$, $\dim (w \cdot J) \cap W=\dim W - \dim \pi_J(W)$.
\end{lem}

\begin{proof}
By applying the fibre dimension theorem to the projection $\pi_J:W \rightarrow (F^\times)^n/J$.
\end{proof}

\begin{lem}\label{final zar closed}
Let $D$ be a full $K$-powered field with a strong finitely generated $K$-powered subfield $D_0$. Suppose $(L,W)$ is a free pair in $D^n$ with $\dim L + \dim W \leqslant n$. Then:

\begin{itemize}
    \item[(a)] There is a Zariski-closed proper subset $W'$ of $W$ such that if $\exp(b) \in \exp(L) \cap (W \setminus W')$, then $\trdeg(\exp(b)/F_0) >0$.
    \item[(b)] There is a Zariski-closed proper subset $W''$ of $W$ such that if $\exp(b) \in \exp(L) \cap (W \setminus W'')$ then $\ldim_\mathbb{Q}(b/D_0)=n$.
\end{itemize}
\end{lem}

These proofs are inspired by Section 5 of the unpublished preprint \cite{Zil11}.

\begin{proof}[Proof of $(a)$] Let $H_0$ denote the maximal $\mathbb{Q}$-affine subspace over $D_0$ contained in $L$; then $H_0$ is a translate of some $\mathbb{Q}$-linear subspace $H$. Let $\pi_H:D^n \twoheadrightarrow D^n/H$ and $\pi_{\exp(H)}: (F^\times)^n \twoheadrightarrow (F^\times)^n/\exp(H)$ denote the projections. Then 
\begin{align*}
    \dim \pi_H(L) + \dim \pi_{\exp(H)}(W) & \leqslant (\dim L - \dim H) + \dim W \\
    & \leqslant n-\dim H.
\end{align*}

By fixing a $\mathbb{Q}$-linear isomorphism between $D^n/H$ and $D^{n-\dim H}$ we may treat the pair $(\pi_H(L),\pi_{\exp(H)}(W))$ as a pair in $D^{n-\dim H}$. 

Assume that $(\pi_H(L),\pi_{\exp(H)}(W))$ satisfies the Lemma, and let $Z$ be the corresponding Zariski-closed proper subset of $\pi_{\exp(H)}(W)$. Let $W':= \pi^{-1}_{\exp(H)}(Z)$, and assume $\exp(b) \in (\exp(L) \cap (W \setminus W'))$. Then $\pi_{\exp(H)}(\exp(b)) \notin Z$, so $\trdeg(\pi_{\exp(H)}(\exp(b))/\pi_{\exp(H)}(D_0)) >0$. Thus, we obtain $\trdeg(\exp(b)/D_0)>0$, as we wanted.

Hence we may assume without loss of generality that $L$ does not contain any positive-dimensional $\mathbb{Q}$-affine subspaces over $D_0$.

Since $D_0$ is finite dimensional as a $\mathbb{Q}$-vector space, $\exp(D_0)$ is a finite rank subgroup of $(F^\times)^n$. By Laurent's Theorem \autocite[Theor\`eme 1]{Lau84} (the multiplicative group case of the Mordell-Lang Conjecture) there is a finite union $M$ of cosets of proper algebraic subgroups of $(F^\times)^n$ such that every point of $W \cap \exp(D_0^n)$ lies in $M$. Since $(L,W)$ is free, $W':=W \cap M$ is a proper Zariski-closed subset of $W$.

Suppose $\exp(b) \in \exp(L) \cap (W \setminus W')$. Then in particular $b \notin D_0^n$, and $\ldim_\mathbb{Q}(b/D_0)>0$.

Let $Q=\mathbb{Q}$-$\affloc(b/D_0)$. Since  we assumed that $L$ does not contain any positive-dimensional $\mathbb{Q}$-affine subspaces over $D_0$, $\dim (L \cap Q) < \dim Q$. Then we have $$\ldim_K(b/D_0) \leqslant \dim (L \cap Q) < \dim Q = \ldim_\mathbb{Q}(b/D_0).$$

Since $D_0 \treq D$, we have that $\delta(b/D_0)\geqslant 0.$ Then \begin{align*}
    \trdeg(\exp(b)/F_0) & \geqslant \ldim_\mathbb{Q}(b/D_0) - \ldim_K(b/D_0)>0. \qedhere
\end{align*} 
\end{proof}

\begin{proof}[Proof of $(b)$] For each $\mathbb{Q}$-linear subspace $H$ of $D^n$, with projections $\pi_H:D^n \rightarrow D^n/H$ and $\pi_{\exp(H)}:(F^\times)^n \rightarrow (F^\times)^n/\exp(H)$, let $W_H$ denote the empty set whenever $\dim (\pi_H(L))+\dim(\pi_{\exp(H)}(W)) > n-\dim H$, and the Zariski-closed subset obtained applying $(a)$ to the pair $(\pi_H(L), \pi_{\exp(H)}(W))$ otherwise (again, we can do this by identifying $D^n/H$ with $D^{n-\dim H}$). Note that for the trivial subspace $H=\{0\}$, $W_{\{0\}}$ is the Zariski-closed set obtained by applying $(a)$ to $(L,W)$ itself.

By Theorem \ref{weak cit} there is a finite set $\mathcal{H}$ of $\mathbb{Q}$-linear subspaces of $D^n$ such that every atypical component of an intersection between $W$ and an algebraic subgroup of $(F^\times)^n$ is contained in a coset of $\exp(H)$ for some $H \in \mathcal{H}$. By applying Lemma \ref{first typical fibre lemma} to each of the projection maps $\pi_{\exp(H)}:W \rightarrow (F^\times)^n/\exp(H)$ for $H \in \mathcal{H}$, and taking the union of all the (finitely many) Zariski-closed sets thus obtained, we find a Zariski-closed proper subset $W_\pi$ of $W$ such that for each $w \in W \setminus W_\pi$ and each $H \in \mathcal{H}$, $\dim (w \cdot \exp(H)) \cap W= \dim W - \dim \pi_{\exp(H)}(W)$.

Let then $$W'':=W_\pi \cup W_{\{0\}} \cup \bigcup_{H \in \mathcal{H}} \pi_{\exp(H)}^{-1}(W_H).$$ Let $b \in D^n$ be a point such that $\exp(b) \in \exp(L) \cap (W \setminus W'')$; we aim to show that $\ldim_\mathbb{Q}(b/D_0)=n$. Let $Q=\mathbb{Q}$-$\affloc(b/D_0)$, and let $S$ be the irreducible component of $W \cap \exp(Q)$ which contains $\exp(b)$. Since $\exp(b) \notin W_{\{0\}}$, $\exp(b)$ is not algebraic over $F_0$, so $\dim S >0$.

Since $S$ is positive dimensional, there is a $\mathbb{Q}$-linear space $H \in \mathcal{H} \cup \{D^n\}$ such that $S \subseteq \exp(b+H)$. By Remark \ref{w is algebraic}, the coset $\exp(b+H)$ is defined over $F_0^\alg$. The preimage of $F_0^\alg$ under $\exp$ is an infinite-dimensional $\mathbb{Q}$-vector subspace of $D$, and hence we may find $c$ in $\exp^{-1}(F_0^\alg)$ such that:

\begin{itemize}
    \item[(i)] $b+H=c+H$,
    \item[(ii)] $\ldim_\mathbb{Q}(c/D_0)=\ldim_\mathbb{Q}(b/D_0)=\dim Q$, and
    \item[(iii)] $\ldim_\mathbb{Q}(c/D_0b)=\dim (Q \cap b+H)$, 
\end{itemize}
from which we obtain that $\ldim_\mathbb{Q}(b/D_0c)=\dim (Q \cap c+H)$.

Moreover, $S$ is a typical component of the intersection $(W \cap (\exp(c+H))) \cap (\exp(Q) \cap \exp(c+H))$ with respect to $\exp(c+H)$: this means that 
\begin{align}
\dim S & = \dim (W \cap \exp(c+H)) + \dim (\exp(Q \cap c+H)) - \dim (\exp(c+H)) \nonumber \\
& = \dim (W \cap \exp(c+H)) + \dim (Q \cap c+H) - \dim H.  \label{eighttwo}
\end{align} 

Note that we are allowing the possibility that $H=D^n$, i.e$.$, that $S$ is typical: we will see that under our assumptions this is the only possible case.

Since $\exp(c)$ is algebraic over $F_0$ and $D_0 \treq D$, $\delta(c/D_0)=0$: this implies that $\delta(b/D_0c)=\delta(b,c/D_0) \geqslant 0$. Then we have 
\begin{align*}
    0 & \leqslant \delta (b/D_0c) \\
    & = \ldim_K(b/D_0c)+\trdeg(\exp(b)/F_0) -\ldim_\mathbb{Q}(b/D_0c) \\
    & \leqslant \dim (L \cap Q \cap (c+H)) + \dim S - \dim ( Q \cap (c+H)).
\end{align*}
Combining this with (\ref{eighttwo}) we find that 
\begin{align*}
    0 & \leqslant \dim (L \cap  Q \cap (c+H)) + \dim (W \cap \exp(c+H)) - \dim H \\
    & \leqslant \dim (L \cap (c+H)) + \dim (W \cap \exp(c+ H))-\dim H.
\end{align*} 
Therefore $$\dim H \leqslant  \dim (L \cap (c+H)) + \dim (W \cap \exp(c+ H)).$$

Since $\exp(b) \notin W_\pi$ and $\exp(c) \cdot \exp(H)=\exp(b) \cdot \exp(H)$, we have $$\dim (\pi_{\exp(H)}(W))=\dim W - \dim (W \cap \exp(c+H)).$$ Hence
\begin{alignat*}{3}
    \dim \pi_{H}(L) + \dim \pi_{\exp(H)}(W)  &= &&(\dim L -   \dim (L \cap c+H))                  & \\
                                               &  && + (\dim W - \dim (W \cap \exp(c+H)))         & \\
                                               &= &&(\dim L +  \dim W )                             & \\
                                               &  && - (\dim (L \cap c+H) + \dim (W \cap \exp(c+H)) & \\
                                               & \leqslant && n - \dim H  & 
\end{alignat*} 

Thus $(\pi_{H}(L),\pi_{\exp(H)}(W))$ satisfies the assumptions the Lemma. However, as the irreducible component $S$ is contained in a translate of $\exp(H)$ by an element that is algebraic over $F_0$, we have that $\trdeg(\pi_{\exp(H)}(\exp(b))/F_0)=0$. As we have taken $b$ so that $\exp(b) \notin W_{H_j}$ for each $H_j \in \mathcal{H}$, by part $(a)$ of the Lemma this is only possible if $H=D^n$. Hence $S$ is a typical component of $W \cap \exp(Q)$, that is, $\dim S = \dim W + \dim Q -n$.

Since $D_0 \treq D$, we have 
\begin{align*}
    0 & \leqslant \delta(b/D_0) \\
    & =\ldim_K(b/D_0) + \trdeg(\exp(b)/F_0) - \ldim_\mathbb{Q}(b/D_0)\\
    & \leqslant \dim (L \cap Q) + \dim S - \dim Q \\
    & \leqslant \dim L + \dim S - \dim Q \\
    & \leqslant n - \dim W + \dim S - \dim Q \\
    & = 0
\end{align*} where the last equality holds by typicality of the intersection. Therefore the chain of inequalities collapses, and $\dim (L \cap Q)=\dim L$: since $(L,W)$ is free, this only holds for $Q=D^n$. Hence, $D^n=\mathbb{Q}$-$\affloc(b/D_0)$, so $\ldim_\mathbb{Q}(b/D_0)=n$, as required.
\end{proof}

This yields the proof that in this setting $K$-powers-closedness implies strong $K$-powers-closedness.

\begin{proof}[Proof of Proposition \ref{pow closed implies spow closed}]
Suppose $D_0$ is a finitely generated strong substructure of the full $K$-powered field $D$; let $(L,W)$ be a free, rotund, $n$-dimensional pair in $D^n$, $a \in D^k$.

Extending $a$ if necessary, we may assume that $\langle D_0, a \rangle_\mathbb{Q} \treq D$. By Lemma \ref{final zar closed} there is a proper Zariski-closed subset $W''$ of $W$ such that if $\exp(b) \in \exp(L) \cap (W \setminus W'')$ then $\ldim_\mathbb{Q}(b/D_0,a)=n$. Since $D$ is $K$-powers-closed, $\exp(L) \cap W$ is Zariski-dense in $W$, and therefore such a point $\exp(b)$ exists. Now apply Proposition \ref{saturation and pow closedness} to deduce $D$ is algebraically saturated in $\mathcal{C}(D_0)$.
\end{proof}

\section{Generic powers}\label{generic powers section}

In this section we explore the consequences of our work for the $K$-powered fields $\mathbb{C}^K$ and $\mathbb{B}^K$ for which we know an appropriate Schanuel statement. For particularly nice $K$, we can prove $\mathbb{C}^K \cong \mathbb{B}^K$. Note that the field $K$ plays different roles. In $\mathbb{E}^K$, or in a general $K$-powered field, it is an abstract field. We define $\mathbb{C}^K$ and $\mathbb{B}^K$ by taking $K$ as a subfield of $\mathbb{C}$ (or $\mathbb{B}$) and taking account of how the exponential interacts with that subfield.

We first give the categoricity and quasiminimality result we have proved. It includes Theorem \ref{third} as the special case in which $D_0$ is the standard base $SB^K$ introduced in Example \ref{standardbase}.

\begin{thm}\label{better third}
	Let $K$ be a countable field of characteristic $0$, $D_0$ a finitely generated $K$-powered field with cyclic kernel. 
	
	Then up to isomorphism there is exactly one $K$-powered field $\mathbb{E}^K(D_0)$ of cardinality continuum which:

\begin{itemize}
	\item[(i)] is a strong extension of $D_0$,
	\item[(ii)] is $K$-powers closed, and
	\item[(iii)] has the countable closure property.
\end{itemize}	

Furthermore, it is quasiminimal.
\end{thm}

\begin{proof}
	By Theorem \ref{characterize} there is only one $K$-powered field which satisfies (i) and (iii), and which is algebraically saturated in the category $\mathcal{C}(D_0)$. By Propositions \ref{saturation and pow closedness} and \ref{pow closed implies spow closed}, algebraic saturation in $\mathcal{C}(D_0)$ is equivalent to $K$-powers closedness.
\end{proof}

\begin{thm}\label{third for B}
    If $K \subseteq \mathbb{B}_\exp$ has finite transcendence degree, then there is a finitely generated partial $K$-powered subfield $D_0 \treq \mathbb{B}^K$ such that $\mathbb{B}^K \cong \mathbb{E}^K(D_0)$.
\end{thm}

\begin{proof}
    By Lemma \ref{fingenstrong} if $K$ has finite transcendence degree then $\mathbb{B}^K$ has a finitely generated strong substructure $D_0 \subseteq \hull{K}_\exp$. Proposition \ref{pclccp} gives the countable closure property. $K$-powers-closedness of $\mathbb{B}^K$ follows from exponential-algebraic closedness of $\mathbb{B}_\exp$ (see \autocite[Section 5]{Zil05} and \autocite[Section 9]{BK18}). Then apply Theorem \ref{better third}.
\end{proof}

We pick out those subfields $K$ where $D_0$ can be taken to be trivial.

\begin{defn}
Let $F$ be $\mathbb{C}_{\exp}$ or $\mathbb{B}_{\exp}$. We say that a subfield $K \subseteq F$ of finite transcendence degree \textit{acts as a subfield of generic powers} if $F^K \cong \mathbb{E}^K$. 

If $\lambda \in F$ is transcendental and such that $\mathbb{Q}(\lambda)$ acts as a subfield of generic powers, we say that $\lambda$ \textit{is a generic power}.
\end{defn}

\begin{exas}\label{taugeneric}
(i) The kernel generator $\tau$ is a generic power. If we consider the predimension $\delta^{\mathbb{Q}(\tau)}$, we obtain:
\begin{align*}
\delta^{\mathbb{Q}(\tau)}(z/\tau) & =\ldim_{\mathbb{Q}(\tau)}(z/\tau) + \trdeg(\exp(z)) - \ldim_{\mathbb{Q}}(z/\tau) \\ 
& \geqslant \trdeg(z/\tau)+ \trdeg (\exp(z)) - \ldim_\mathbb{Q}(z/\tau) \\
& \geqslant \trdeg(z,\exp(z)/\tau) - \ldim_{\mathbb{Q}}(z/\tau) \\
& = \delta^{\exp}(z/\tau) \\
& \geqslant 0
\end{align*}
which implies that $\hull{\varnothing}_{\mathbb{Q}(\tau)}=SB^{\mathbb{Q}(\tau)}$ and hence $\mathbb{B}^{\mathbb{Q}(\tau)} \cong \mathbb{E}^{\mathbb{Q}(\tau)}$.

(ii) A similar argument shows that the number $\pi=\frac{\tau}{2i}$ in $\mathbb{B}_\exp$ is a generic power.

(iii) Not all transcendental $\lambda \in \mathbb{B}_\exp$ are generic powers. Let $\lambda_1=\log \sqrt{2}, \lambda_2=\log \sqrt{3}, \lambda=\frac{\lambda_1}{\lambda_2}$. 

Then:
\begin{align*}
\delta^{\mathbb{Q}(\lambda)} (\lambda_1,\lambda_2/\tau) &= \ldim_{\mathbb{Q}(\lambda)}(\lambda_1,\lambda_2/\tau) + \trdeg(\exp(\lambda_1), \exp(\lambda_2)) - \ldim_{\mathbb{Q}}(\lambda_1,\lambda_2/\tau)\\
&= 1 + 0 -2\\
&=-1.
\end{align*}
Therefore $\mathbb{B}^{\mathbb{Q}(\lambda)} \ncong \mathbb{E}^{\mathbb{Q}(\lambda)}$.
\end{exas}

The question of isomorphism between $\mathbb{C}_\exp$ and $\mathbb{B}_\exp$ is considered out of reach as it requires Schanuel's Conjecture. However, for ``sufficiently generic'' tuples of complex numbers, we can prove the corresponding statement for powers.

\begin{thm}\label{theorem on generic powers}
    Let $K=\mathbb{Q}(\lambda_1,\dots,\lambda_n)$ be the field of rational functions, and choose embeddings of $K$ into $\mathbb{C}_\exp$ and into $\mathbb{B}_\exp$ such that the $\lambda_i$ are exponentially-algebraically independent.

    Then $K$ acts as a field of generic powers, so we have $\mathbb{C}^K \cong \mathbb{E}^K \cong \mathbb{B}^K$.
\end{thm}

\begin{proof}
    By \autocite[Theorem 1.3]{BKW} we have $SB^K \treq \mathbb{C}^K$ and $SB^K \treq \mathbb{B}^K$. Hence, by Theorem \ref{third for B}, $\mathbb{B}^K \cong \mathbb{E}^K$. $\mathbb{C}^K$ is $K$-powers-closed by \autocite[Corollary 8.10]{Gal22}, and it has the countable closure property by Corollary \ref{pclccp}. Therefore, by Theorem \ref{better third} $\mathbb{C}^K \cong \mathbb{E}^K$.
\end{proof}

Theorem \ref{fourth} follows from Theorem \ref{theorem on generic powers}.

\section{Quasiminimality}\label{quasiminimality section}

In this last section we focus on the category $\mathcal{C}^\tr(D_0)$, for $D_0$ a countable full $K$-powered field: we show that in this category algebraic saturation is equivalent to a notion of \textit{generic strong $K$-powers-closedness}, which in turn follows from $K$-powers-closedness as we show using the Weak CIT again, although in a different way. As a consequence, we obtain that every $K$-powers-closed field with the countable closure property (in particular $\mathbb{C}^K$ for any countable subfield $K$ of $\mathbb{C}$, and thus $\mathbb{C}^\mathbb{C}$) is quasiminimal. The method of the proof and the terminology are similar to the ones used in \autocite[Section 11]{BK18} to reduce Conjecture \ref{qm conj} to exponential-algebraic closedness, but the technical details are different.

We note that the usage of the word ``generic'' in \textit{generic $K$-powers-closedness} is not related to the usage in the notion of \textit{generic power} discussed in the previous section. 

\begin{defn}\label{gspc definition}
Suppose $D_0$ is a full countable $K$-powered field and $D$ is a purely powers-transcendental extension of $D_0$.

Let $(L,W)$ be a free, rotund, $n$-dimensional pair in $D^n$ that is defined over some $a \in D^k$ such that $\ldim_{\mathbb{Q}}(a/D_0)=k$, and such that for all $(c,\exp(c))$ generic in $L \times W$ over $D_0$ we have that $\Loc(c,a/D_0)$ is free and strongly rotund.

$D$ is \textit{generically strongly $K$-powers-closed over $D_0$} if for every $(L,W)$ and $a$ as above, there is $(b,\exp(b)) \in L \times W$ such that $\ldim_{\mathbb{Q}}(b/D_0a)=n$.
\end{defn}

We have an analogue of Proposition \ref{saturation and pow closedness}, replacing $\mathcal{C}(D_0)$ by $\mathcal{C}^\tr(D_0)$.

\begin{prop}\label{ppt-saturation}
Let $D_0$ be a full countable $K$-powered field. Recall that $\mathcal{C}^\tr(D_0)$ denotes the full subcategory of $\mathcal{C}(D_0)$ whose objects are purely powers-transcendental extensions of $D_0$.

An object in $\mathcal{C}^\tr(D_0)$ is algebraically saturated in $\mathcal{C}^\tr(D_0)$ if and only if it is full and generically strongly $K$-powers-closed over $D_0$.
\end{prop}

\begin{proof}
    $(\Rightarrow)$ Let $D$ be an object of $\mathcal{C}^\tr(D_0)$ that is algebraically saturated in $\mathcal{C}^\tr(D_0)$. It is straightforward to show that $D$ is full, so we show that it is generically strongly $K$-powers-closed.
    
    Let $a \in D^k$ with $\ldim_{Q}(a/D_0)=k$, and let $(L,W)$ be a pair as in Definition \ref{gspc definition}. Replacing $a$ by a basis for $\hull{D_0, a}_K^D$ if necessary, we may assume that $D_1:=\langle D_0, a \rangle_{\mathbb{Q}} \treq D$. Add to $D_1$ a point $b$ such that $(b,\exp(b))$ is generic in $L \times W$ over $D_1$; by Proposition \ref{classify}, this generates a strong, powers-algebraic extension $D_2$ of $D_1$, which is purely powers-transcendental over $D_0$ and hence is an object of $\mathcal{C}^\tr(D_0)$. By algebraic saturation of $D$, $D_2$ embeds into $D$ over $D_1$. Hence there is a point $b'$, the image of $b$ under the embedding, such that $(b',\exp(b')) \in L \times W$, and (by genericity) $\ldim_\mathbb{Q}(b'/D_1)=n$. Hence $D$ is generically strongly $K$-powers-closed.

    $(\Leftarrow)$ Let $D$ be an object of $\mathcal{C}^\tr(D_0)$ that is full and generically strongly $K$-powers-closed over $D_0$. Let $D_1=\langle D_0, a\rangle_\mathbb{Q}$ be a finitely generated, purely powers-transcendental extension of $D_0$ that is strong in $D$, and let $D_2$ be a finitely generated, strong, powers-algebraic extension of $D_1$ that is purely powers-transcendental over $D_0$. By Proposition \ref{goodbasesprop}, we may find a good basis $b$ for $D_2$ over $D_1$, so $D_2=\langle D_1, b \rangle_\mathbb{Q}$.

    We prove that $D_2$ strongly embeds into $D$ over $D_1$ by induction on $n:=\dim_{\mathbb{Q}}(D_2/D_1)$. Let $(L,W)=\Loc(b/D_1)$ and $(L_1,W_1)=\Loc(b,a/D_0)$. 
    
    If $(L,W)$ is free, then since the extension $D_1 \treq D_2$ is strong and powers-algebraic, we have that $(L,W)$ is rotund and $n$-dimensional. Moreover, since $D_0 \treqcl D_2$, $(L_1, W_1)$ is strongly rotund; by genericity of $b$ over $D_1$ then $(L,W)$ satisfies the assumptions in Definition \ref{gspc definition}, so by generic strong $K$-powers-closedness of $D$ there is a point $b' \in D^n$ such that $(b', \exp(b')) \in L \times W$ and $\ldim_{\mathbb{Q}}(b'/D_1)=n$. Since $D$ strongly extends $D_1$, it must be the case that $\ldim_K(b'/D_1)+\trdeg(\exp(b')/F_1) \geqslant n$. As $\dim L + \dim W =n$, this implies that $(L,W)=\Loc(b'/D_1)$. Since $b$ is a good basis, the extension $\langle D_1, b' \rangle_\mathbb{Q}$ of $D_1$ is then isomorphic to $D_2$ over $D_1$, and it is contained in $D$, so $D_2$ embeds in $D$ over $D_1$.

    If $(L,W)$ is not free, then without loss of generality we may assume that either $L$ is contained in a $\mathbb{Q}$-affine subspace over $D_1$ of the form $z_i=c$ for some $c \in V_1$, or $W$ is contained in an algebraic subvariety of the form $w_i=d$ for some $d \in F_1^\alg$. In the first case we replace $D_1$ by $D_1':=\hull{D_1,c}_K^D$ and consider $D_2$ as a powers-algebraic strong extension of $D_1$ of lower dimension (hence we may apply induction); the second case is treated analogously.
\end{proof}

We also have the analogue of Proposition \ref{pow closed implies spow closed}. Some of the ideas in the proof are similar, but overall the argument is different.

\begin{prop}\label{pcimpgpc}
Let $D$ be a full, purely powers-transcendental extension of a countable full $K$-powered field $D_0$.

If $D$ is $K$-powers-closed, then it is generically strongly $K$-powers-closed.
\end{prop}

Before proving this proposition, we prove two lemmas about atypical intersections, the second a uniform version of the first.

\begin{lem}\label{nonunif}
Let $(L,W)$ be a free, strongly rotund pair in $D^n$. Suppose $W$ is defined over an algebraically closed subfield $F_0$ of $F$.

For every $\mathbb{Q}$-linear proper subspace $Q \leq D^n$ there is a Zariski-closed proper subset $W' \subseteq W$, defined over $F_0$, such that if $w \in W \setminus W'$, then $$\dim ( W \cap w \cdot \exp(Q)) < \dim W + \dim L -n -\dim (Q \cap L)+ \dim Q. $$
\end{lem}

This is similar to an atypical intersection statement, but we do not require the intersection $W \cap w \cdot \exp(Q)$ to be typical - it is allowed to be atypical within a certain margin, determined  by the interplay of $Q$ with $L$.

\begin{proof}
Let $\pi_Q:D^n \twoheadrightarrow D^n/Q$ and $\pi_{\exp(Q)}:(F^\times)^n \twoheadrightarrow (F^\times)^n/\exp(Q)$ denote the projections. By strong rotundity of the pair we have $$\dim \pi_Q(L)+ \dim \pi_{\exp(Q)} (W) > n-\dim Q.$$

By the fibre dimension theorem, there is a Zariski-closed proper subset $W'_\pi$ of $\pi_{\exp(Q)}(W)$, defined over $F_0$, such that for all $w \cdot \exp(Q) \in \pi_{\exp(Q)}(W) \setminus W'_\pi$, $$\dim \pi_{\exp(Q)}^{-1}(w \cdot \exp(Q)) = \dim W - \dim \pi_{\exp(Q)}(W).$$ Let $W'=\pi_{\exp(Q)}^{-1}(W'_\pi)$. We know $\pi_{\exp(Q)}^{-1}(w \cdot \exp(Q))$ is by definition equal to $W \cap w \cdot \exp(Q)$ for all $w \in W$. On the other hand, since $\pi_Q(L)$ is a linear projection of a linear space, the dimension of the fibres of the points in its image does not depend on the choice of point, and it is always equal to $\dim (L \cap Q)=\dim L - \dim \pi_Q(L)$.

Combining the fibre equalities, we obtain that for $w \in W \setminus W'$, 
\begin{align*}
    \dim (L \cap Q) + \dim (W \cap w \cdot \exp(Q)) =  \dim& L - \dim \pi_Q(L) +\\
     & \dim W - \dim \pi_{\exp(Q)}(W)
\end{align*}  and thus using the inequality obtained by strong rotundity 
\begin{align*}
    \dim (L \cap Q) + \dim (W \cap w \cdot \exp(Q)) &<  \dim L + \dim W-n + \dim Q
\end{align*}
from which we obtain the statement.
\end{proof}

Using the finiteness given by the weak CIT we show that the set $W'$ can in fact be chosen uniformly in $Q$.

\begin{lem}\label{uniform}
Let $(L,W)$ be a free, strongly rotund pair. Suppose $W$ is defined over an algebraically closed subfield $F_0$ of $F$.

There is a Zariski-closed proper subset $W'$ of $W$, defined over $F_0$, such that for every $\mathbb{Q}$-linear proper subspace $Q$ of $D^n$ and every $w \in W\setminus W'$, $$\dim (W \cap w \cdot \exp(Q)) <\dim W + \dim L -n -\dim (Q \cap L)+\dim Q. $$
\end{lem}

\begin{proof}
By Theorem \ref{weak cit}, there is a finite list $\mathcal{H}=\{H_1,\dots,H_l\}$ of $\mathbb{Q}$-linear spaces of $D^n$ such that every positive dimensional atypical component of an intersection between $W$ and a coset of an algebraic subgroup of $(F^\times)^n$ is contained in a coset of $\exp(H_j)$ for some $H_j \in \mathcal{H}$.

For each $\mathbb{Q}$-linear subspace $H \in \mathcal{H}$, let $W'_H$ denote the Zariski-closed proper subset defined over $F_0$ given by Lemma \ref{nonunif}. Let $W':=\bigcap_{H \in \mathcal{H}}W'_H$. This is an intersection of sets defined over $F_0$, so it is also defined over $F_0$.

Suppose $w \in W \setminus W'$, let $Q \leq D^n$ be a proper $\mathbb{Q}$-linear space, and let $X$ be the irreducible component of $W \cap w \cdot \exp(Q)$ containing $w$. 

Assume $X$ is a typical component of the intersection, so that it satisfies $\dim X = \dim W + \dim Q -n$. Then, since $(L,W)$ is free and $Q$ is a proper subspace, $\dim L - \dim L \cap Q>0$ and therefore 
\begin{align*}
    \dim X &< \dim W + \dim L -n - \dim L\cap Q + \dim Q
\end{align*}
as required. Thus, we assume that $X$ is atypical.

Then there is a $\mathbb{Q}$-linear space $H \in \mathcal{H}$ such that 
\begin{align}\label{seven}
    \dim X &= \dim W \cap w \cdot \exp(H) + \dim Q \cap H - \dim H.
\end{align} 

Since $w \notin W'_{H_j}$ for any $H_j \in \mathcal{H}$, and thus in particular $w \notin W'_H$, we must have
\begin{align}\label{eight}
    \dim W \cap w \cdot \exp(H) &< \dim W + \dim L - n - \dim L \cap H + \dim  H
\end{align} The combination of (\ref{seven}) and (\ref{eight}) gives 
\begin{align}
    \dim X & < \dim W + \dim L - n - \dim L \cap H + \dim Q \cap H.\label{nine}
\end{align}

Since $\dim L \cap H  \geqslant \dim L \cap Q \cap H$, 
\begin{align}
    \dim Q \cap H - \dim L \cap H & \leqslant \dim Q \cap H - \dim L \cap (Q \cap H) \nonumber \\
    & = \dim (L+ (Q \cap H)) - \dim L \nonumber \\
    & \leqslant \dim (L + Q) - \dim L \nonumber \\
    &= \dim Q - \dim L \cap Q. \label{thirteen}
\end{align}

So combining (\ref{nine}) and (\ref{thirteen}), 
\begin{align*}
    \dim X < \dim W + \dim L - n - \dim L \cap Q + \dim Q 
\end{align*}
 in this case as well.
\end{proof}

Now we can prove that $K$-powers-closedness implies generic strong $K$-powers-closedness.

\begin{proof}[Proof of Proposition  \ref{pcimpgpc}]
Suppose $D$ is a full $K$-powered field which is $K$-powers closed and is a purely powers-transcendental extension of a countable $K$-powered subfield $D_0$. Let $(L,W)$ be a free, rotund, $n$-dimensional pair in $D^n$ defined over $D_0a$ for some $a \in D^k$ that is as in the definition of generic strong $K$-powers-closedness. Then for $(c,\exp(c))$ generic in $L \times W$ over $D$ we have that $(L_1,W_1):=\Loc(c,a,\exp(c,a)/D_0)$ is free and strongly rotund in $D^{n+k}$. Let $W_1'$ be the Zariski-closed proper subset of $W_1$ defined over $F_0$ given by Proposition \ref{uniform}. Let $W_1^\circ:=W_1 \setminus W_1'$. 

Consider the set $$W^\circ:=\left\{ w \in W \mid (w,\exp(a)) \in W^\circ_1 \right\}.$$ Since $\exp(a)$ is generic over $F_0$ in the projection of $W_1$ to the last $k$ coordinates, $W^\circ$ is a Zariski-open dense subset of $W$. By $K$-powers-closedness of $D$, there is a point $(b,\exp(b)) \in L \times W^\circ$; note that then $\exp(b,a) \in W_1^\circ$ and $(b,a) \in L_1$. We will prove that $\ldim_\mathbb{Q}(b/D_0a)=n$. 

Let $Q:=\affloc_{\mathbb{Q}}(b/D_0a)$ and $Q_1:=\affloc_{\mathbb{Q}}(b,a/D_0)$. Extending $a$ if necessary, we assume that $D_0 a \treq D$. Therefore we have
\begin{alignat}{3}
    0 & \leqslant \delta(b/D_0a)                          &&                                                         & \nonumber\\
      & = \ldim_K(b/D_0a)+\trdeg(\exp(b)/F_0         && \exp(a)) - \ldim_\mathbb{Q}(b/D_0a)                     & \nonumber\\
      & = \ldim_K(b,a/D_0) - \ldim_K(a/D_0)          && + \trdeg(\exp(b,a)/D_0) - \trdeg(\exp(a)/D_0)           &\nonumber \\
      &                                              &&- \ldim_\mathbb{Q}(b,a/D_0) + \ldim_\mathbb{Q}(a/D_0)    &\nonumber \\
      & \leqslant \dim(L_1\cap Q_1) + \dim (W_1 \cap \exp(&&Q_1))  - \dim Q_1                                        & \nonumber\\
      &                                              && - (\ldim_K(a/D_0)+\trdeg(\exp(a)/F_0)-k)                &\nonumber \\
      & = \dim (L_1 \cap Q_1) + \dim (W_1 \cap \exp( &&Q_1))   - \dim Q_1 - \delta(a/D_0).                      & \label{fourteen}
\end{alignat}
From the fibre dimension theorem we get that $\ldim_K(a/D_0)=\dim L_1 - \dim L$ and $\trdeg(\exp(a)/F_0)=\dim W_1-\dim W$, and we also have that $\dim L + \dim W =n$, thus 
\begin{align}
    \delta(a/D_0)&=\dim L_1 + \dim W_1 -n -k \label{fifteen}
\end{align} Therefore, combining (\ref{fourteen}) and (\ref{fifteen}) we obtain $$0 \leqslant \dim (L_1 \cap Q_1) + \dim (W_1 \cap \exp(Q_1)) - \dim Q_1 - (\dim L_1 + \dim W_1 - n -k)$$ so $$\dim (W_1 \cap \exp(Q_1)) \geqslant \dim L_1 - \dim (L_1 \cap Q_1) + \dim W_1 -n -k + \dim Q_1.$$

Recall that $\exp(Q_1)$ is a coset of some algebraic subgroup by $\exp(b,a)$: since $\exp(b,a) \in W_1^\circ$, this inequality together with Proposition \ref{uniform} implies that $Q_1$ is not a proper subspace, so $Q_1=D^{n+k}$. Then $Q=D^n$, and $\ldim_\mathbb{Q}(b/D_0a)=n$, as we wanted.
\end{proof}

Finally we put everything together to prove our main quasiminimality theorems.

\begin{proof}[Proof of Theorem \ref{second}]
Let $K$ be a countable subfield of $\mathbb{C}$, and $D_0=\pcl_K^\mathbb{C}(\varnothing)$. Then $\mathbb{C}^K$ is $K$-powers-closed by \autocite[Corollary 8.10]{Gal22}, so by Propositions \ref{ppt-saturation} and \ref{pcimpgpc} it is algebraically saturated in $\mathcal{C}^\tr(D_0)$. It has the countable closure property by Proposition \ref{pclccp}, and it has cardinality continuum, so by Theorem \ref{characterize} it is isomorphic to $\mathbb{E}^{K,\tr}(D_0)$, and hence is quasiminimal.
\end{proof}

As noted in the Introduction, Theorem \ref{first} follows directly from Theorem \ref{second}.

\printbibliography

\textsc{Francesco Gallinaro, Mathematisches Institut, Albert-Ludwigs-Universit\"at Freiburg, Ernst-Zermelo-Str. 1, 79100 Freiburg, Germany.}

E-mail address: \texttt{francesco.gallinaro@mathematik.uni-freiburg.de}.

\vskip 10mm

\textsc{Jonathan Kirby, School of Mathematics, University of East Anglia, Norwich Research Park, NR4 7TJ, Norwich, UK.}

E-mail address: \texttt{jonathan.kirby@uea.ac.uk}.

\end{document}